\documentclass[a4paper,11pt]{article}
\usepackage[]{times}
\usepackage{fullpage}
\usepackage{graphicx}
\usepackage{float}
\usepackage{subfigure}
\usepackage{amsmath}
\usepackage{fancyhdr}
\usepackage[colorlinks]{hyperref}
\usepackage{xcolor}
\date{}

\newcommand{\prov}{{\sc Proof}.\hspace*{0mm} }

\newcommand{\QED}{$\rule{2mm}{2mm}$}

\newcommand{\natu}{{\sf I \! N}}
\newtheorem{theorem}{Theorem}[section]
\newtheorem{lemma}[theorem]{Lemma}
\newtheorem{e-proposition}[theorem]{Proposition}

\newtheorem{e-definition}[theorem]{Definition\rm}
\newtheorem{remark}{\it Remark\/}

\title{Optimal enforcement of natural conditions on smooth boundaries \\ with piecewise polynomials upon fitted straight-edged meshes}
\author{
     Vitoriano Ruas$^{1}$\thanks{Sorbonne Universit\'e, Campus Pierre et Marie Curie, 4 place jussieu, Couloir 55-65, 4\`eme \'etage, 75005 Paris, France.}
		\\[1mm]
  {\small $^{1}$ Institut Jean Le Rond d'Alembert, CNRS UMR 7190, Sorbonne Universit\'e, Paris, France.}\\[1mm]
  {\small e-mail: {\it vitoriano.ruas@upmc.fr}}}

\begin{document}
\maketitle

\begin{abstract}
A few decades ago some possible remedies to an inaccurate enforcement of Neumann or Robin conditions prescribed on the boundary of a smooth domain, owing to its approximation by the union of straight-edged triangles or tetrahedra in a fitted mesh, were addressed. These studies,  due to authors such as Barrett \& Elliott \cite{BarrettElliott}, advocated the use of isoparametric finite elements with a single curvilinear edge or face fitting the true boundary at two or three vertexes and also at additional points on those curves or curved surfaces, so as to define non polynomial piecewise solutions with the optimal order of approximation. In this work we adopt a different approach, whose main feature is the use of a fitted mesh consisting only of straight-edged elements, upon which the solution is a polynomial. Lost optimal orders of convergence, by virtue of the domain's approximation by a polytope, are recovered by means of the addition of terms to the bilinear form. These account for natural boundary conditions of the same type as the true ones, though to be prescribed on the approximating boundary instead. The new technique is applied to the case of triangular Lagrange finite elements, for which we give a rigorous reliability study in the solution of reaction-diffusion equations. Numerical experimentation is supplied in support of the theoretical results.
        

\end{abstract}

\section{Preamble}

To begin with, the author endeavors to situate this contribution of his to the numerical enforcement of boundary conditions with enhanced accuracy. More specifically, the focus here is on the solution by the finite element method or related techniques of problems with natural conditions prescribed on a smooth boundary. A paper's outline is also provided at the end of this section. 

\subsection{A short historical account} 

Accurate non-affine simplicial finite-element solution of second order boundary value problems posed in curved domains was a topic of extensive investigation in the early stages of the mathematical studies of this popular method. Most of the contributions to the topic, such as those by Nitsche \cite{Nitsche}, Bramble et al. \cite{BrambleDT} and Scott \cite{Scott}, attempted to use polynomial shape- and test-functions in the most frequent case of simplicial meshes. However those approaches were confined to two-dimensional problems with Dirichlet boundary conditions. In parallel, by that time, the isoparametric finite-element technique based on curvilinear elements and non affine transformations thereof onto a straight-edged master element, was by far the most widespread among practitioners (cf. \cite{Zienkiewicz}) to tackle such a problem. This is because it allowed to handle both two- and three-dimensional problems with either essential or natural boundary conditions. Later authors such as Ciarlet \& Raviart \cite{CiarletRaviart}, Nedoma \cite{Nedoma} and Barrett \& Elliott \cite{BarrettElliott}, to quote just a few, provided a rigorous theoretical background to justify this method's reliability in different contexts. \\
In spite of such outstanding and celebrated work, more recently this topic has gained interest with the introduction of techniques bypassing non affine transformations. This is because the latter require a judicious and not always obvious choice of numerical integration procedures, among other drawbacks. Furthermore isoparametric finite elements cannot be viewed as a universal technique, since they have been mostly developed for second order boundary value problems in standard formulation and, some decades later, also in mixed form \cite{FBRT}. 
Among more recent contributions to the field, the Shifted Boundary Method can be mentioned \cite{ACScovazzi} as a promising one. We could also quote a parameter free evolution of Nitsche's \cite{Nitsche} and Bramble et al.'s \cite{BrambleDT} approaches proposed in \cite{DupontGS}. It should be noted that recent contributions to the field seem to have focused mostly on the case where sole Dirichlet conditions are prescribed on curvilinear boundaries. In other words, a lot of effort to study finite-element techniques of optimal order to handle Neumann or Robin boundary conditions in the case of a smooth domain was done for isoparametric approximations a long time ago in the aforementioned articles \cite{Nedoma} and \cite{BarrettElliott}, among others quoted therein. 

\subsection{Methodological background}     
    
In the framework of essential conditions, the author and collaborators addressed a new approach to handle more accurately DOFs (degrees of freedom) prescribed on the boundary of a smooth $N$-dimensional domain for $N=2,3$ in a series of publications finalized within the last ten years. In contrast to the isoparametric method, in this approach both the shape- and test-function spaces consist only of piecewise polynomials and the computational domain is the polytope formed by the union of straight-edged $N$-simplexes of a fitting mesh. But its key point is the use of a trial space - or a manifold in the case of non zero DOFs -, different from the test space, in the sense that, prescribed DOFs on the boundary are enforced at their exact locations for the former, while such values are enforced at shifted locations on the boundary of the approximating polytope for the latter. We refer to \cite{Arxiv2017}, \cite{JAMP}\cite{Maugin}, \cite{ZAMM}, \cite{IMAJNA}, \cite{M2AN} and references therein for mote details on this methodology, as applied to different kinds of problems and formulations. \\
\indent 
It is noteworthy that similar principles have been exploited in parallel by other authors, mostly in connection with the finite volume method (see e.g. \cite{Clain}). More recently they were also applied in \cite{DG2D} to the DG (discontinuous Galerkin) method in the two-dimensional case. The reliability of the resulting method restricted to Dirichlet boundary conditions was formally established in \cite{AAMS}, by adapting to the DG environment the mathematical analysis of this type of technique provided in \cite{ZAMM} for the finite element method. We also observe that in these works such an approach is referred to as the ROD-method, where the acronym stands for \textit{reconstruction of off-site data}. In contrast, in our own contributions such as \cite{ZAMM}, we call it a Petrov-Galerkin formulation, thereby emphasizing that it is characterized by the use of different spaces of trial- and test-functions. As a matter of fact, the methodology to be studied here is inspired by the approach adopted in \cite{ZAMM}, though adapted to the case of Neumann or Robin boundary conditions. However, as seen below, in this case there is no need to work with different trial- and test-spaces, even though boundary conditions are mimicked at their exact locations in a very similar way to the case of Dirichlet conditions. For this reason it is also  appropriate to use the term \textit{off-site enforcement of boundary conditions} to characterize the method addressed in this work.\\ Incidentally, it should be emphasized that the technique advocated in this work can be applied as such to enhance the quality of Discontinuous Galerkin approximations of the same type of problems, at the price of minor modifications inherent to this method.  \\
\indent
As already pointed out, as far as the finite-element approximation of order greater than one of Neumann or Robin problems 
is concerned, mostly in the eighties, authors like Barrett and Elliott \cite{BarrettElliott}, \v{Z}eni\v{s}ek \cite{Zenisek} and Cermak \cite{Cermak} advocated the use of elements with a curved edge in the two-dimensional case. In doing so, the computational domain comes closer to the actual domain, so that the expected order of approximation remains unchanged as compared to the case of a polygonal domain. To the best of our knowledge, since then no rigorous studies on new approaches to handle this problem appeared in the literature. For this reason we believe that this work brings about an effective alternative for the accurate treatment of natural conditions prescribed on curvilinear boundaries, as we attempt to show hereafter. 

\subsection{Paper's outline}    

An outline of the article is as follows. In Section 2 we present the motivation of our new technique, by describing it in main lines in the framework of a simple model second order elliptic PDE. In Section 3, after giving some general notations, definitions and assumptions to be used throughout the article, we describe our method as applied to a model linear reaction-diffusion equation with Neumann or Robin boundary conditions in a smooth domain. In Section 4 we formally establish underlying stability and convergence results in the natural norm for both convex and non convex domains. Numerical experimentation on a comparative basis with other solution techniques is reported in Section 5; as a by-product, we validate the theoretical reliability results. We conclude in Section 6 with some comments and perspectives for future work.     
        \section{Introduction}
 
\hspace{4mm} In order to clarify the principle the method roduced in this work is based upon, nothing is better than the illustration of the issue to be resolved in a practical case: Assume that we wish to solve by the finite element method the following reaction-diffusion equation with constant strictly positive diffusion and reaction coefficients whose ratio is $c_0$, in a smooth convex two-dimensional domain $\Omega$ with boundary $\Gamma$, namely, \\
Given a function $f_0$ defined in $\Omega$, whose regularity matches in some sense that of its definition domain, we wish to find $u_0$ fulfilling\\  
\begin{equation}
\label{eq0}
- \Delta u_0 + c_0 u_0 = f_0,
\end{equation}
and satisfying homogeneous Neumann boundary conditions on the whole $\Gamma$. \\
The most used finite element method to solve equation \eqref{eq0} is based on continuous functions, which are polynomials of degree not greater than a certain integer $k >0$ in each straight-edged triangle $T$, viewed as a closed set in a fitted mesh ${\mathcal T}_h$ of $\Omega$. The latter expression means that all the vertexes of the polygonal domain $\Omega_h$ lie on $\Gamma$, where $\Omega_h$ is the interior of the set $\displaystyle \cup_{T \in {\mathcal T}_h}$ with boundary $\Gamma_h$ .\\
First of all we represent by ${\mathcal A} \cdot {\mathcal B}$ the inner product of two tensors ${\mathcal A}$ and ${\mathcal B}$ of order $n$ for $n \in \natu$. Then for all functions $w$ and $v$ in the Sobolev space $H^k(\Omega)$ (cf. \cite{Adams}) for $k=1$, we define a bilinear form $a_{0h}$ and a linear form $L_{0h}$ by
\begin{equation}
\label{blf}
a_{0h}(w,v) := \int_{\Omega_h} (\nabla w \cdot \nabla v \;+ \; c_0 wv) \mbox{ and } L_{0h}(v) := \int_{\Omega_h} f_0 v. 
\end{equation}
The variational problem whose solution is the underlying finite element approximation $u_{0h}^k$ of $u_0$ writes as follows:
\begin{equation}
\label{vp0h}
\mbox{Find } u^k_{0h} \in V^k_h \mbox{ such that } a_{0h}(u^k_{0h},v) = L_{0h}(v) \; \forall v \in V_h^k,
\end{equation}
where $V^k_h$ is the space of continuous functions defined in $\Omega_h$, whose restriction to every $T \in {\mathcal T}_h$ is a polynomial of degree less than or equal to $k > 1$. \\
It is well known that problem \eqref{vp0h} has a unique solution. Moreover, an upper bound for the approximation error in the standard norm $\| \cdot \|_{1,h}$ of $H^1(\Omega_h)$ is given by (cf. \cite{COAM})
\begin{equation}
\label{erb}
\| u_0 - u^k_{0h} \|_{1,h} \leq \displaystyle \frac{1}{\alpha_0} \left[ A_0 \displaystyle \inf_{w \in V^k_h} \| u_0 - w \|_{1,h} + \displaystyle 
\sup_{ v \in V^k_h \setminus \{0\}} \frac{a_{0h}(u_0,v)-L_{0h}(v)}{\| v \|_{1,h}} \right], 
\end{equation}
where $\alpha_0 = \min[1,c_0]$ and $A_0 = \max[1,c_0]$. \\
Denoting by $h_T$ the maximum edge length of $T \in {\mathcal T}_h$, in the usual setting that a regular family of meshes is in use, defining $h:= \displaystyle \max_{T \in {\mathcal T}_h} h_T$, the \textit{inf}-term in \eqref{erb} can be estimated by an $O(h^k)$-term, provided $u \in H^{k+1}(\Omega)$. \\
The case of the \textit{sup}-term instead deserves more attention, owing to the mismatch of $\Omega$ and $\Omega_h$. Actually in the remainder of this work we focus on the numerator of this term, commonly called the \textit{variational residual} of the approximation method. In order to estimate it properly we make the very reasonable assumption that no triangle in ${\mathcal T}_h$ has more than two vertexes on $\Gamma$. Let ${\mathcal S}_h$ be the subset of ${\mathcal T}_h$ consisting of triangles having exactly two vertexes on $\Gamma$. For every $T \in {\mathcal S}_h$, $e_T$ represents the edge of $T$ contained in $\Gamma_h$. The portion of $\Gamma$ comprised between the two end-points of $e_T$ is denoted by $\Gamma_T$. Further, ${\bf n}$ and ${\bf n}_h$ being the unit outer normal vectors to $\Gamma$ and $\Gamma_h$, we denote by $\partial \cdot/\partial n$ and $\partial \cdot/\partial n_h$ the first order outer normal derivatives to $\Gamma$ and $\Gamma_h$. \\
Now, using integration by parts, we have
\begin{equation}
\label{vr1}
a_{0h}(u_0,v)-L_{0h}(v) = \displaystyle \int_{\Omega_h} (- \Delta u_0 + c_0 u_0 - f_0)v + \oint_{\Gamma_h} \frac{\partial u_0}{\partial n_h}v, 
\end{equation}
that is,
\begin{equation}
\label{vr2}
a_{0h}(u_0,v)-L_{0h}(v) = \oint_{\Gamma_h} \frac{\partial u_0}{\partial n_h}v \; \forall v \in V^k_h,  
\end{equation}
or yet,
\begin{equation}
\label{vr3}
a_{0h}(u_0,v)-L_{0h}(v) = \displaystyle \sum_{T \in {\mathcal S}_h} \int_{e_T} \frac{\partial u_0}{\partial n_h}v \; \forall v \in V^k_h.  
\end{equation}
By assumption we know that $\partial u_0/\partial n \equiv 0$ on $\Gamma$, but in general $\partial u_0/\partial n_h \neq 0$. Thus the true order of the method will result from the estimation of the magnitude of the latter normal derivative. Let us go into it.\\
To begin with, we denote by $\Delta_T$ the closed set with the smallest area delimited by $e_T$ and $\Gamma$ $\forall T \in {\mathcal S}_h$. Now, referring to Figure 1 and assuming that $h$ is sufficiently small, to every $M \in e_T$ it is possible to associate a single point $P \in \Gamma_T$ as the intersection with $\Gamma$ of the perpendicular to $e_T$ passing through $M$. Clearly enough, the length of the segment $\overline{MP}$ is bounded above by $C_{\Gamma} h_T^2$, where $C_{\Gamma}$ is a constant independent of $T$ (see e.g. \cite{LewNegri}).\\
\begin{figure}[h]
\label{fig1}
\centerline{\includegraphics[width=3.0in]{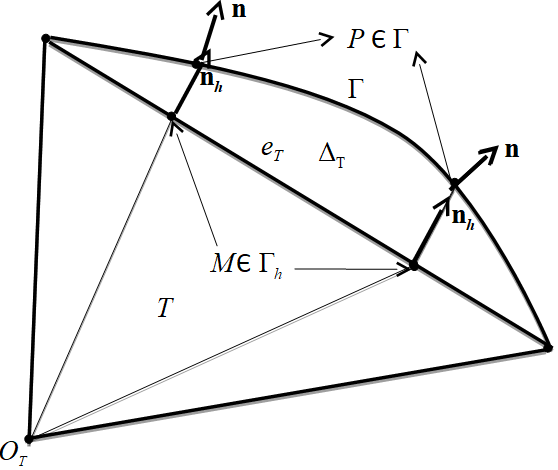}}
\vspace*{8pt}
\caption{Two points $P \in \Gamma$ associated with two points $M \in e_T$ for $T \in {\mathcal S}_h$}
\end{figure}
Now, denoting by $\bf{t}$ the tangent vector along $\Gamma$ oriented as a rotation of ${\bf n}$ by $\pi/2$ in the counterclockwise sense, and by $\partial \cdot/ \partial \tau$ the tangential derivative along $\Gamma$ in the direction of $\bf{t}$, we have \\ 
\begin{equation}
\label{dudnh1} 
\displaystyle \frac{\partial u_0}{\partial n_h}(M) = \frac{\partial u_0}{\partial n_P}(M) \mbox{ cos } \theta_T - \frac{\partial u_0}{\partial \tau_P}(M) \mbox{ sin } \theta_T,
\end{equation}
where $\theta_T$ is the angle $({\bf n}_h,{\bf n}_P)$ restricted to $e_T$ and $\partial (\cdot)/\partial n_Q$ and $\partial (\cdot)/\partial \tau_Q$ denote respectively the partial derivatives in the fixed directions ${\bf n}_Q$ and ${\bf t}_Q$ of the normal and tangential vectors to $\Omega$ at a given point $Q \in \Gamma$. \\
On the other hand, since by assumption $[\partial u_0/ \partial n](P) =0$, the first partial derivative on the right hand side of \eqref{dudnh1} can be rewritten as
\begin{equation}
\label{dudnh2}
\displaystyle \frac{\partial u_0}{\partial n_P}(M) = - \displaystyle \int_M^P \frac{\partial^2 u_0}{\partial n_P \partial n_h}. 
\end{equation} 
It is clear that sin $\theta_T$ is bounded above by $C_{\theta} h_T$, where $C_{\theta}$ is a constant independent of $T$ (cf. \cite{ZAMM}). Plugging \eqref{dudnh2} into \eqref{dudnh1} we obtain
\begin{equation}
\label{dudnh3} 
\displaystyle \left| \frac{\partial u_0}{\partial n_h}(M) \right| \leq \displaystyle \left| \int_M^P \frac{\partial^2 u_0}{\partial n_P \partial n_h} \right| + C_{\theta} \left|\frac{\partial u_0}{\partial \tau_P}(M) \right| h_T.
\end{equation}
 Now, in order to simplify the estimation, let us assume that $u_0$ belongs to the Sobolev space $H^{3+\delta}(\Omega)$ for a small $\delta >0$. This implies that each component of the Hessian of $u_0$ denoted by $H(u_0)$ belongs to the space $L^{\infty}(\Omega)$ with norm $\| \cdot \|_{\infty}$ and that moreover there exists a constant $C_{\delta}$ such that 
$$\| H(u_0) \|_{\infty} \leq C_{\delta} \| u_0 \|_{3+\delta},$$ 
where the notation $\| \cdot \|_{s}$ stands for the standard norm of $H^{s}(\Omega)$ for any strictly positive real number $s$.
On the other hand, taking into account \eqref{dudnh3} and recalling \eqref{vr3} we have 
\begin{equation}
\label{vr4}
a_{0h}(u_0,v)-L_{0h}(v) \leq 
\max[C_{\theta},C_{\Gamma}] \displaystyle \sum_{T \in {\mathcal S}_h} h_T \left[ \int_{e_T} \left|\frac{\partial u_0}{\partial \tau_P} v \right| +  h_T \| H(u_0) \|_{\infty} \int_{e_T} |v| \right],      
\end{equation}
which after straightforward calculations leads to
\begin{equation}
\label{vr5}
a_{0h}(u_0,v)-L_{0h}(v) \leq 
\max[C_{\theta},C_{\Gamma}] \left( h \| \nabla u_0 \|_{\infty}  +  h^2 \| H(u_0) \|_{\infty} \right)\int_{\Gamma_h} |v|.      
\end{equation} 
By the Trace theorem we have $\int_{\Gamma_h} |v| \leq C(\Gamma_h) \| v \|_{1,h}$, where the constant $C(\Gamma_h)$ is assumed to be independent of $h$ (cf. Remark 1 hereafter). Hence it follows from \eqref{vr5} that, in principle, the sup-term in \eqref{erb} is bounded above by an $O(h)$-term. Eventually, under particular conditions such as $\partial u_0/\partial \tau \equiv 0$ on $\Gamma$, a higher order estimation of this term might hold, but in any case the best we can hope for is an upper bound with an $O(h^2)$-term. This is certainly a serious limitation for the standard finite-element solution method \eqref{vp0h}, at least for $k > 2$. \\
As a matter of fact, the above premises themselves suggest a variant of \eqref{vp0h} that bypasses the aforementioned order limitation. The idea behind it is the famous idiom:\\ 
"If the mountain doesn't come to you then you must go to the mountain".\\
Transposing it to our specific case, if we wish to annihilate the variational residual, why don't we simply subtract from the bilinear form $a_{0h}(w,v)$ the problematic term $$\displaystyle -\int_{e_T} \left[\mbox{ sin } \theta_T \frac{\partial w}{\partial \tau_P} + \mbox{ cos } \theta_T \int_M^P \frac{\partial^2 w}{\partial n_P \partial n_h} \right] v$$ $\forall P \in \Gamma_T$ and $\forall T \in {\mathcal S}_h$? This does work indeed, as we endeavor to show in the remainder of this article in a rigorous and more general framework at a time.      



\section{More accurate handling of natural conditions on smooth boundaries}

In this section we describe the finite element formulation that we advocate for the solution of a model second order elliptic equation 
in a smooth domain $\Omega \subset \Re^2$, with inhomogeneous Neumann or Robin boundary conditions. In contrast to the previous section, henceforth $\Omega$ is no longer assumed to be convex. \\

\subsection{Notations, definitions and assumptions on the mesh}

Before going into the description of our method, we specify the notations to be used in the sequel, in addition to those already introduced in Section 2. \\
We begin with notations related to well known Sobolev function spaces (see e.g. \cite{Adams}).\\
$D$ being a bounded subset of $\Re^2$, $\parallel \cdot \parallel_{s,D}$ and $| \cdot |_{s,D}$ denote the standard norm and semi-norm of Sobolev space $H^{s}(D)$, for $s \in \Re^+$ with $H^0(D)=L^2(D)$. We further denote by $\parallel \cdot \parallel_{m,p,D}$ and $| \cdot |_{m,p,D}$ the usual norm and semi-norm of $W^{m,p}(D)$ for $m \in \natu^{*}$ and $p \in [1,\infty] \setminus \{2\}$ with $W^{0,p}(D)=L^p(D)$, and eventually for $W^{m,2}(D)=H^m(D)$ as well. In case $D=\Omega$ we omit this subscript. \\
As for additional notations related to the mesh ${\mathcal T}_h$,  the polygonal domain $\Omega_h$ 
we extend $\| \cdot \|_{1,h}$ to the standard norm of either $H^j(\Omega_h)$ for $j > 1$ or $L^2(\Omega_h)$ for $j=0$. Moreover, in order to accommodate non convex domains, we introduce the set $\tilde{\Omega}_h := \Omega \cup \Omega_h$. We also represent by $\tilde{T}$ the set $T \cup \Delta_T$ $\forall T \in {\mathcal S}_h$ and define $\tilde{T} := T \; \forall T \in {\mathcal T}_h \setminus {\mathcal S}_h$.
In this work the notation $|\cdot|$ is employed to represent both the area of a subset of $\Re^2$ with a non empty interior and the length of a segment or a curve.\\     
$P_k(D)$ is defined to be the space of polynomials of degree less than or equal to $k$ restricted to a subset $D \in \Re^2$. 
\\

To conclude we list below some additional assumptions on the meshes.\\
Besides being fitted, unlike previous work on the subject, the family of meshes ${\mathcal T}_h$ under consideration denoted by ${\mathcal P}$ is not assumed to be a quasi-uniform family of triangulations of $\Omega$. Indeed, in this article we only require ${\mathcal P}$ to be a    regular family of meshes. Recalling the notation $h_T$ introduced in Section 2 and denoting the shortest edge of $T$ by $l_T$, this means that there exists a strictly positive constant $\mu$ such that $l_T \geq \mu h_T$, $\forall T \in {\mathcal T}_h$ and $\forall {\mathcal T}_h \in {\mathcal P}$. \\
Referring to Figure 1, throughout this article we assume that, for all meshes in ${\mathcal P}$, the indexing parameter $h:= \displaystyle \max_{T \in {\mathcal T}_h} h_T$ is sufficiently small for the intersection $P$ with $\Gamma_T$ of the perpendicular to $e_T$ passing through $M \in e_T$ to be unique $\forall M \in e_T$ and $\forall T \in {\mathcal S}_h$.\\ 
Since $h$ is supposedly small we may also assume that ${\mathcal P}$ is such that there exist two mesh-independent constants $s_m$ and $c_m$ fulfilling 
\begin{equation}
\label{Thetam}
\left\{
\begin{array}{l}
\max_{T \in {\mathcal S}_h} \max_{M \in e_T} | \theta_T(M)| \;  \leq \; s_m h \implies 
\max_{T \in {\mathcal S}_h} \max_{M \in e_T} | sin \theta_T(M)| \; \leq \; s_m h \\  
\mbox{and } \min_{T \in {\mathcal S}_h} \min_{M \in e_T} | cos \theta_T(M)| \; \geq \; c_m > 0 \; \forall {\mathcal T}_h \in {\mathcal P}.
\end{array}
\right.
\end{equation} 
\begin{remark} In complement to this preliminary material we make an important remark on practical consequences of the fact that $h$ is assumed to be small. More particularly, it is about constants appearing in celebrated inequalities for Sobolev spaces defined in $\Omega_h$. Since these constants  depend on $\Omega_h$, in principle they depend on $h$. However, even though we do not explicitly prove it here for the sake of conciseness, we shall take for granted without notice, that they have upper bounds independent of $h$, wherever they appear in the sequel. Nevertheless we note that, in the same manner as in the Supplementary Material of \cite{M2AN}, it is possible to prove the existence of such upper bounds, since $\Omega_h$ is just a small perturbation of $\Omega$. A significant example thereof is the constant $C_{tr}$ of the trace inequality \eqref{trace} given hereafter. \QED
\end{remark}   

\subsection{The model problem}
Next we introduce the model problem that we chose for the study of our methodology.   
In order to focus on essential aspects, we consider the same reaction-diffusion equation with inhomogeneous Robin boundary conditions as in the work by Barrett and Elliott \cite{BarrettElliott}. This choice brings about simplifications, in that many results used by these authors allow us to bypass some rather cumbersome estimations in our analysis. 

\begin{remark} In the numerical experimentation section we apply our method to reaction-advection-diffusion equations with a variable advective velocity. In doing so, we show that this does not really makes any difference, as far as error estimates are concerned. \rule{2mm}{2mm}
\end{remark}

Let us assume that $\Omega$ is a smooth domain and that, for an integer $k >1$,    
$f$ and $g$ are given in $H^{k-1}(\Omega)$ and $H^{k-1/2}(\Gamma)$. We also assume that $d$ is a function defined in $\bar{\Omega}$ 
fulfilling $d \geq d_0 > 0$, $d_0 \in \Re$, and also that $r$ and $\eta$ are functions defined in $\Omega$ and $\Gamma$ such that $r \geq r_0 \geq 0$ and $\eta \geq \eta_0 \geq 0$, where $r_0$ and $\eta_0$ are real numbers satisfying $r_0 + \eta_0 >0$. Moreover, since our method is basically designed to deal with quite smooth solutions, we shall assume that $\Omega$, $d$, $r$ and $\eta$ enjoy all regularity properties necessary for the results in the theoretical analysis performed hereafter to hold true. \\
Now the problem to solve is   
\begin{equation}
\label{eq}
\left\{
\begin{array}{l}
\mbox{Find } u \in H^{k+1}(\Omega) \mbox{ such that} \\
- \nabla \cdot ( d \nabla u) + r u = f \; \mbox{ in } \Omega \\
\mbox{and } \displaystyle d \frac{\partial u}{\partial n} + \eta u = g \mbox{ on } \Gamma.
\end{array}
\right. 
\end{equation}
Recalling the domain $\tilde{\Omega}_h$ we define $\tilde{\Omega}$ to be the union of $\tilde{\Omega}_h$ for all triangulations ${\mathcal T}_h \in {\mathcal P}$. We observe that $\tilde{\Omega} = \Omega$ if $\Omega$ is convex. However, in the general case $\tilde{\Omega}_h$ is not a subset of $\Omega$, so that $\Omega$ is a strict subset of $\tilde{\Omega}$. Therefore, if $\Omega$ is not convex, a priori it is necessary to extend the data $f$, $d$ and $r$ in $\tilde{\Omega} \setminus \Omega$ to functions $\tilde{f} \in H^{k-1}(\tilde{\Omega})$, $\tilde{d} \in C^k(\bar{\tilde{\Omega}})$ and $\tilde{r} \in C^{k-1}(\bar{\tilde{\Omega}})$ in order to enable the proper definition of an approximate problem for the domain $\Omega_h$ in all cases. We assume that there exists a strictly positive constant $\tilde{d}_0$ and a non negative constant $\tilde{r}_0$ such that $\tilde{d} \geq \tilde{d}_0$ and $\tilde{r} \geq \tilde{r}_0$ all over $\tilde{\Omega}$. For convenience we also consider that $u$ is extended to the whole $\tilde{\Omega}$ by a function $\tilde{u} \in H^{k+1}(\tilde{\Omega})$, for example, in the way advocated in \cite{Stein}, among others. Notice, however, that $\tilde{u}$ is by no means expected to be a solution of a reaction-diffusion equation like \eqref{eq} in $\tilde{\Omega}$ with $d$, $r$ and $f$ replaced with $\tilde{d}$, $\tilde{r}$ and $\tilde{f}$. However, as seen hereafter, the knowledge of $\tilde{f}$ outside $\Omega$ is not indispensable. Thus we will assume that $\tilde{f}:= -\nabla \cdot (\tilde{d}\nabla \tilde{u}) + \tilde{r}\tilde{u}$ in $\tilde{\Omega}\setminus \Omega$, even if $\tilde{u}$ is unknown in $\tilde{\Omega} \setminus \Omega$.     

\subsection{Finite-element analog}

Recalling the finite-dimensional space $V^k_h$ introduced in Section 2, in order to set up the finite-element analog advocated in this work to solve \eqref{eq}, we first extend every $v \in V^k_h$ to $\Omega \setminus \Omega_h$ as follows: $\forall T \in {\mathcal S}_h$ the expression of $v$ in $T$ in terms of the space variables is applied as such to all points in $\Delta_T$ not belonging to $T$. We represent by $\tilde{V}^k_h$ the space of continuous functions $\tilde{v}$ defined in $\tilde{\Omega}_h$ such that $\forall T \in {\mathcal T}_h$ the restriction of $\tilde{v}$ to $\tilde{T}$ belongs to $P_k(\tilde{T})$.\\
Now, referring to Figure 1, with every function $\varphi$ defined on $\Gamma$ we associate a function $\bar{\varphi}$ uniquely defined on $\Gamma_h$ by $\bar{\varphi}(M) = \varphi(P)$ $\forall M \in e_T$ and $\forall T \in {\mathcal S}_h$, where $P$ is the point of $\Gamma$ associated with $M \in e_T$ in the way specified in Subsection 3.1.\\    
Now, recalling that $\theta_T$ is the angle between ${\bf n}_h$ and ${\bf n}$ restricted to $\Delta_T$ for $T \in {\mathcal S}_h$, let us define for all functions $w \in H^3(\tilde{\Omega}_h) + \tilde{V}_h^k$ and $v \in V_h^k$ a bilinear form $\bar{a}_{h}$ and 
a linear form $\bar{L}_{h}$ as follows.
\begin{equation}
\label{barah}
\left\{
\begin{array}{l}
\bar{a}_h(w,v) := \int_{\Omega_h} (\tilde{d} \nabla w \cdot \nabla v \;+ \; \tilde{r} wv) \\   
+ \displaystyle \sum_{T \in {\mathcal S}_h} \int_{e_T} \left\{\mbox{ cos } \theta_T \left[ \bar{\eta} \bar{w} \; + \; \displaystyle \int_M^P \frac{\partial}{\partial n_h}\left(\tilde{d} \frac{\partial w}{\partial n_P}\right)\right] \; + \; \tilde{d} \mbox{ sin } \theta_T \frac{\partial w}{\partial \tau_P}  \right\}v
\end{array}
\right.
\end{equation}
\begin{equation}
\label{barLh}
\bar{L}_h(v) := \int_{\Omega_h} \tilde{f} v \; + \; \displaystyle \sum_{T \in {\mathcal S}_h} \int_{e_T} \mbox{ cos } \theta_T \bar{g} v.
\end{equation}
A natural way to define our finite element approximation of \eqref{eq} is to 
solve the following problem
\begin{equation}
\label{barvph}
\mbox{Find } \bar{u}^k_h \in \tilde{V}^k_h \mbox{ such that } \bar{a}_{h}(\bar{u}^k_h,v) = \bar{L}_h(v) \; \forall v \in V^k_h. 
\end{equation}
However, in general, evaluating both $\bar{a}_h$ and $\bar{L}_h$ exactly may be unpractical if not unfeasible. Therefore, we pose instead a similar variational formulation in which suitable approximation techniques such as numerical integration or interpolation are employed in $T$ for $T \in {\mathcal T}_h$ or in $e_T$ for $T \in {\mathcal S}_h$. Clearly enough, the deviation from the exact solution brought about by this additional error source must be compatible with the expected order of convergence of the formulation \eqref{barvph}. \\
As for $\bar{a}_h$, in every $T \in {\mathcal T}_h$ we use a $J_k$-point Gauss quadrature formula with points $Q_j \in T$ and strictly positive weights $\lambda_j$ for $j=1,\ldots J_k$ for non polynomial integrands in $T$ and an $I_k$-point Gaussian quadrature formula for integrals along $e_T$ with points $M_i$ and strictly positive weights $\omega_i$ for $i=1,\ldots I_k$, whose sum equals one. We will be more specific about both formulae in due course. For the moment we just approximate integrals of a given continuous function $\phi$ in $T \in {\mathcal T}_h$ by ${\mathcal J}_T^k(\phi)$, where 
\begin{equation}
\label{JTk} 
{\mathcal J}_T^k(\phi) := |T| \displaystyle \sum_{j=1}^{J_k} \lambda_j \phi(Q_j) 
\end{equation}
In order to do without the explicit knowledge outside of $\Omega$ of the extensions $\tilde{d}$ and $\tilde{r}$ of the data $d$ and $r$, we assume that for all $T \in {\mathcal S}_h$ all the quadrature points $Q_j$ lie inside $\Omega$, in case this domain is not convex.\\
The integrals along $e_T$ for $T \in {\mathcal S}_h$ in turn are approximated by ${\mathcal I}_T^k(\phi)$, where           
\begin{equation}
\label{ITk} 
{\mathcal I}_T^k(\phi) := |e_T| \displaystyle \sum_{i=1}^{I_k} \omega_i \phi(M_i). 
\end{equation}  
Now for every $w \in H^3(\tilde{\Omega}_h) + \tilde{V}_h^k$ and $v \in V_h^k$ we set
\begin{equation}
\label{ah}
a_h(w,v) := b_h(w,v) + c_h(w,v), 
\end{equation}
where $\forall (w;v) \in [H^3(\tilde{\Omega}_h) + \tilde{V}^k_h] \times V^k_h$
\begin{equation}
\label{bh}
b_h(w,v) := \; \displaystyle \sum_{T \in {\mathcal T}_h} {\mathcal J}_T^k \left( \tilde{d} \nabla w \cdot \nabla v \; + \; \tilde{r} w v \right) \; 
+ \; \displaystyle \sum_{T \in {\mathcal S}_h} {\mathcal I}_T^k \left(\mbox{cos } \theta_T \bar{\eta} w v \right)  
\end{equation}
and
\begin{equation}
\label{ch}
c_h(w,v) := \; \displaystyle \sum_{T \in {\mathcal S}_h} {\mathcal I}_T^k  \left[ \left\{\mbox{cos } \theta_T \left[ \bar{\eta} (\bar{w}-w) \; + \; \displaystyle \int_M^P \frac{\partial}{\partial n_h}\left(\tilde{d} \frac{\partial w}{\partial n_P}\right)\right] \; + \; \tilde{d} \mbox{ sin } \theta_T \frac{\partial w}{\partial \tau_P} \right\} v \right]. 
\end{equation}
In contrast, as far as $\bar{L}_h$ is concerned, we define its approximation $L_h$ by
\begin{equation}
\label{Lh}
L_{h}(v) := \int_{\Omega_h} \rho_{k-1}(\tilde{f}) v \; + \; \displaystyle \sum_{T \in {\mathcal S}_h} \int_{e_T} \sigma_{k-1}(\mbox{ cos } \theta_T \bar{g})_{|e_T} v \; \forall v \in V^k_h.
\end{equation}
where for an integer $l$ fulfilling $k-1 \geq l \geq 1$, $[\rho_l(\phi)]_{|T} \in P_{l}(T)$ is a standard Lagrange interpolate of a continuous function $\phi$ at $(l+1)(l+2)/2$ distinct points $R_j$ of $T$ such that no $R_j$ lies in $T \setminus \Omega$ whenever this set is non empty; $[\sigma_{l} (\phi)]_{|e_T} \in P_l(e_T)$ in turn is the standard Lagrange interpolate of $\phi$ at $l+1$ equally spaced points $S_i$ of $e_T$ including its end-points.\\   
With the above definitions, we set the approximate problem to solve as
\begin{equation}
\label{vph}
\mbox{Find } u^k_{h} \in \tilde{V}^k_h \mbox{ such that } a_{h}(u^k_{h},v) = L_{h}(v) \; \forall v \in V_h^k,
\end{equation}
Notice that both $V^k_h$ and $\tilde{V}^k_h$ are finite-dimensional spaces with the same dimension, Hence both $a_h$ and $L_h$ are continuous forms over $\tilde{V}^k_h \times V^k_h$ and $V^k_h$ equipped with any norm. For example, we have
\begin{equation}
\label{continuityah}
a_h(w,v) \leq A_h \| w \|_{1,h} \| v \|_{1,h} \; \forall (w;v) \in \tilde{V}^k_h \times V^k_h
\end{equation}
and
\begin{equation}
\label{continuityLh}
L_h(v) \leq C_h \| v \|_{1,h} \; \forall v \in V^k_h,
\end{equation}
for suitable constants $A_h$ and  $C_h$.\\
According to \cite{Babuska1} and \cite{COAM}, it follows that problem \eqref{vph} is well posed if the bilinear form $a_h$ satisfies an $inf-sup$ condition with a constant $\alpha_h > 0$ over $\tilde{V}^k_h$ and $V^k_h$ equipped with any norm. In the case under study the natural norm is $\| \cdot \|_{1,h}$, for which this condition writes  
\begin{equation}
\label{infsuph}
\forall w \in \tilde{V}^k_h \displaystyle \sup_{v \in V^k_h \setminus \{0\}} \frac{a_h(w,v)}{\| v \|_{1,h}} \geq \alpha_h \| w \|_{1,h}.
\end{equation}      
Notice that conditions similar to \eqref{continuityah} and \eqref{continuityLh} hold for $\bar{a}_h$ and $\bar{L}_h$. In the next section  we prove condition \eqref{infsuph} together with its counterpart for $\bar{a}_h$, which implies that problem \eqref{vph} is well-posed, as much as \eqref{barvph}.\\

\section{Reliability study of the approximate problem}

The aim of this section is to establish that, besides being well-posed, problem \eqref{vph} satisfies all the conditions required for 
generating approximations of the solution $u$ of \eqref{eq} - or of a natural extension $\tilde{u}$ of $u$ to the whole $\Omega_h$ in case $\Omega_h \setminus \Omega \neq \emptyset$ -, with optimal order $k$ in the norm $\| \cdot \|_{1,h}$ for $k > 1$.   \\
Since our method is designed for a computational domain $\Omega_h$ equal to the union of straight-edged triangles, throughout this section we consider only the particular case where the mappings from the triangles in ${\mathcal T}_h$ onto the master element $\hat{T}$ with vertexes $(0;0)$, $(1;0)$, $(0,1)$ are affine. This allows us to use some results known to hold for isoparametric mappings from triangles with a curved edge onto $\hat{T}$ given in several works on the same subject, such as \cite{CiarletRaviart}, \cite{Nedoma}, \cite{Cermak} and \cite{Zenisek}. \\
Hereafter we use the following notations in connection with the master triangle $\hat{T}$: \\
For every $T \in {\mathcal T}_h$, ${\mathcal B}_T$ is the invertible affine mapping such that $\hat{T} = {\mathcal B}_T(T)$. Further, for every function $\phi$ defined in $T$, $\hat{\phi}$ is the function defined in $\hat{T}$ by $\hat{\phi} = \phi \circ [{\mathcal B}_T]^{-1}$. Finally, for a generic point $Q \in T$, the point of $\hat{T}$ defined by ${\mathcal B}_T(Q)$ is denoted by $\hat{Q}$.\\
We basically need to establish both a uniform \textit{inf-sup} condition - i.e. the uniform weak coercivity - and the uniform continuity of $a_h$, but incidentally, the coercivity of $a_h$ will also be addressed. Furthermore, as an auxiliary material, the same study will be carried out for the bilinear form $\bar{a}_h$  

\subsection{Uniform coercivity}

To begin with, we establish the uniform coercivity of the bilinear form $\bar{a}_h$ over the space $V^k_h$ equipped with the $\| \cdot \|_{1,h}$-norm. With this aim we split it into the sum of two bilinear forms $\bar{b}_h$ and $\bar{c}_h$ given by
\begin{equation}
\label{barbh}
\bar{b}_h(w,v) := \int_{\Omega_h} (\tilde{d} \nabla w \cdot \nabla v \;+ \; \tilde{r} wv) + \displaystyle \sum_{T \in {\mathcal S}_h} \int_{e_T} \mbox{cos } \theta_T \bar{\eta} w v \; \forall (w;v) \in [H^1(\Omega_h)]^2
\end{equation}
and  
\begin{equation}
\label{barch}
\left\{
\begin{array}{l}
\bar{c}_h(w,v):= \; \displaystyle \sum_{T \in {\mathcal S}_h} \int_{e_T} {\mathcal F}_T(w,v) \; \forall (w;v) \in [H^3(\tilde{\Omega}_h) + \tilde{V}^k_h] \times V^k_h \\
\mbox{with} \\
{\mathcal F}_T(w,v) = \displaystyle \left\{\mbox{cos } \theta_T \left[\bar{\eta} (\bar{w}-  w) \; + \; \displaystyle \int_M^P \frac{\partial}{\partial n_h}\left(\tilde{d} \frac{\partial w}{\partial n_P}\right) \right] \; + \; \tilde{d} \mbox{ sin } \theta_T \frac{\partial w}{\partial \tau_P} \right\} v.
\end{array}
\right.
\end{equation}
First we prove 
\begin{e-proposition}
\label{coercivebarbh}
The bilinear form $\bar{b}_h$ is uniformly coercive over $V^k_h$ equipped with the $\| \cdot \|_{1,h}$-norm, in the sense that there exists a mesh-independent constant $\bar{\beta} > 0$ such that 
\begin{equation}
\label{coercivbarbh}
\bar{b}_h(v,v) \geq \bar{\beta} \| v \|_{1,h}^2.
\end{equation}
\end{e-proposition}

\prov
Clearly enough we have 
\begin{equation}
\label{coercivbarbh1}
\bar{b}_h(v,v) \geq \bar{\beta}_1 \| v \|_{1,h}^2 \; \forall v \in V_h^k \mbox{ with } \beta_1 = \min[d_0,r_0].
\end{equation} 
Moreover, recalling \eqref{Thetam}, according to Lemma 2.1 of \cite{Cermak}, for some mesh-independent constant $c >0 $ it holds
\begin{equation}
\label{coercivbarbh2}
\bar{b}_h(v,v) \geq \bar{\beta}_2 \| v \|_{1,h}^2 \; \forall v \in V_h^k \mbox{ with } \beta_2 = \min[d_0,c c_m \eta_0].
\end{equation}
\eqref{coercivbarbh} is thus a consequence of \eqref{coercivbarbh1} and \eqref{coercivbarbh2} with $\bar{\beta} = \max[\bar{\beta}_1,\bar{\beta}_2] >0$. \QED\\

As for $\bar{c}_h$ we have

\begin{e-proposition}
\label{boundbarch}
$\bar{c}_h$ is an $O(h^{1/2})$-bounded bilinear form over $\tilde{V}^k_h \times V^k_h$ equipped with the $\| \cdot \|_{1,h}$-norm, in the sense that there exists a mesh-independent constant $\bar{C}_c > 0$ such that 
\begin{equation}
\label{estimbarch}
\bar{c}_h(w,v) \leq \bar{C}_c h^{1/2} \| w \|_{1,h} \| v \|_{1,h} \; \forall w \in \tilde{V}^k_h \mbox{ and } \forall v \in V^k_h.  
\end{equation}
\end{e-proposition}
\prov
Let $w \in \tilde{V}_h^k$ and $v \in V^k_h$. From \eqref{barch} we have
\begin{equation}
\label{barchwv1}
\left\{
\begin{array}{l}
\bar{c}_h(w,v) \leq \displaystyle \sum_{T \in {\mathcal S}_h}\bar{c}_T(w,v) \\
\mbox{where} \\
\bar{c}_T(w,v)= \displaystyle \int_{e_T} \left\{ \mbox{cos } \theta_T \left[\bar{\eta} \left(  \bar{w}- w \right) \; + \; \displaystyle \int_M^P \frac{\partial}{\partial n_h}\left(\tilde{d} \frac{\partial w}{\partial n_P}\right) \right] \; 
+ \; \tilde{d} \mbox{ sin } \theta_T \displaystyle \frac{\partial w}{\partial \tau_P} \right\}v .
\end{array}
\right.
\end{equation}
Using the identity $[\bar{w}-w](M) = \int_M^P \partial w/\partial n_h$ for $P \in \Gamma$ with $\overline{MP} \perp e_T$, together with the obvious inequality $\| \int_M^P \tilde{\varphi} \|_{0,e_T} \leq \sqrt{C_{\Gamma}} h_T \| \tilde{\varphi} \|_{0,\Delta_T}$ for every square integrable function $\tilde{\varphi}$ in $\tilde{T}$, we easily obtain the following upper bound for the right hand side of the above expression of $\bar{c}_T(w,v)$ with a constant $C_1$ independent of $T$.
\begin{equation}
\label{barchwv2}
\left\{
\begin{array}{l}
\bar{c}_T(w,v) \leq C_{1} \displaystyle \left[ h_T \| \eta \|_{0,\infty,\Gamma} \displaystyle \left\| \frac{\partial w}{\partial n_h}\right\|_{0,\Delta_T} \; + \; h_T  \| \tilde{d} \|_{1,\infty,\tilde{\Omega}} \; \displaystyle \left( \left\| \frac{\partial w}{\partial n_P} \right\|_{0,\Delta_T}^2 \right. \right.\\
\left. \left. + \displaystyle \left\| \frac{\partial^2 w}{\partial n_h \partial n_P} \right\|_{0,\Delta_T}^2 \right)^{1/2} + \;\| \tilde{d} \|_{0,\infty,\tilde{\Omega}} h_T  \left\| \displaystyle \frac{\partial w}{\partial \tau_P} \right\|_{0,e_T}\right] \| v \|_{0,e_T}.
\end{array}
\right.
\end{equation}
Denoting by $\partial \varphi/\partial \tau_h$ the first order partial derivative of a function $\varphi$ along $e_T$, it is clear that 
\begin{equation}
\label{partialtauP}
\partial (\cdot)/\partial \tau_P = \mbox{cos } \theta_T \partial (\cdot)/\partial \tau_h - \mbox{sin } \theta_T \partial (\cdot)/\partial n_h.
\end{equation}
Hence, referring to inverse inequalities for $\tilde{T}$ given in \cite{ZAMM}, it is also easy to figure out that there exists two constants $C_{\tau}$ and $C_{\tau}^{'}$ depending on $\Gamma$ but not on $T$ such that 
\begin{equation}
\label{Ctau}
h_T \displaystyle \left\| \frac{\partial w}{\partial \tau_P} \right\|_{0,e_T} \leq C^{'}_{\tau} \left( h_T^{3/2} \displaystyle \left\|\frac{\partial w}{\partial \tau_h} \right\|_{0,\infty, \tilde{T}} + h_T^{5/2} \displaystyle \left\| \frac{\partial w}{\partial n_h} \right\|_{0,\infty, \tilde{T}} \right) \leq C_{\tau} 
h_T^{1/2} \| w \|_{1,\tilde{T}}.
\end{equation} 
Similarly, we have 
\begin{equation}
\label{partialnP}
\partial (\cdot)/\partial n_P = \mbox{cos } \theta_T \partial (\cdot)/\partial n_h + \mbox{sin } \theta_T \partial (\cdot)/\partial \tau_h.
\end{equation} 
Since $\theta_T$ does not vary in the direction of ${\bf n}_h$, we may also write
\begin{equation}
\label{partial2nP} 
\partial^2 (\cdot)/\partial n_h \partial n_P = \mbox{cos } \theta_T \partial^2 (\cdot)/\partial n_h^2 + \mbox{sin } \theta_T \partial (\cdot)/\partial n_h \partial \tau_h.
\end{equation}
It follows that 
\begin{equation}
\label{prelim}
\displaystyle  \left( \left\| \frac{\partial w}{\partial n_P} \right\|_{0,\Delta_T}^2 + \left\| \frac{\partial^2 w}{\partial n_h \partial n_P} \right\|_{0,\Delta_T} \right)^{1/2} \leq \displaystyle  \left\| \frac{\partial w}{\partial n_h} \right\|_{1,\Delta_T} + s_m h_T \displaystyle  \left\| \frac{\partial w}{\partial \tau_h} \right\|_{1,\Delta_T}.
\end{equation}       
Moreover, resorting to well known results (see e.g. \cite{StrangFix}), since $w_{|T}$ is a polynomial, we know that for a constant $C_n^{'}$ independent of $T$ it holds 
\begin{equation}
\label{GStrangGFix}
 \displaystyle \left\| \frac{\partial w}{\partial n_h} \right\|_{1,\Delta_T} + s_m h_T \displaystyle  \left\| \frac{\partial w}{\partial \tau_h} \right\|_{1,\Delta_T} \leq C_n^{'} h_T^{1/2} \displaystyle \left(\left\| \frac{\partial w}{\partial n_h} \right\|_{1,T} + s_m h_T \displaystyle  \left\| \frac{\partial w}{\partial \tau_h} \right\|_{1,T} \right).
\end{equation}       
Hence, from \eqref{prelim} and \eqref{GStrangGFix}, we easily infer the existence of a constant $C_{n}$ independent of $T$ such that
\begin{equation}
\label{Cn}
h_T \displaystyle  \left( \left\| \frac{\partial w}{\partial n_P} \right\|_{0,\Delta_T}^2 + \left\| \frac{\partial^2 w}{\partial n_h \partial n_P} \right\|_{0,\Delta_T} \right)^{1/2} \leq C_n h_T^{1/2} \| w \|_{1,T}    
\end{equation} 
Collecting \eqref{Ctau}, \eqref{Cn} and recalling another inverse inequality for $\tilde{T}$ given in \cite{ZAMM}, we combine both equations with \eqref{barchwv2} to quite readily obtain for a constant $C_2$ independent of $T$
\begin{equation}
\label{barchwv3}
\bar{c}_T(w,v) \leq C_{2} \displaystyle \left( h_T \| \eta \|_{0,\infty,\Gamma} \; + \; 
h_T^{1/2} \| \tilde{d} \|_{1,\infty,\tilde{\Omega}} \right) \| w \|_{1,T} \| v \|_{0,e_T}.
\end{equation}
Summing up over ${\mathcal S}_h$ and applying the Cauchy-Schwarz inequality to the resulting right hand side, we immediately come up with 
a constant $\tilde{C}(\eta,d)$ depending on $\| \eta \|_{0,\infty,\Gamma}$ and $\| \tilde{d} \|_{1,\infty,\tilde{\Omega}}$ but not on $T$, 
such that
\begin{equation}
\label{barchwv5} 
\bar{c}_h(w,v) \leq \tilde{C}(\eta,\tilde{d}) h^{1/2}  \| v \|_{0,\Gamma_h}  \| w \|_{1,h}.
\end{equation}
Finally, since we can take for granted the existence of a constant $C_{tr}$ independent of $h$ and $v$ such that 
\begin{equation}
\label{trace}
\| v \|_{0,\Gamma_h} \leq \| v \|_{1/2,\Gamma_h} \leq C_{tr} \| v \|_{1,h},
\end{equation} 
\eqref{barchwv5} leads to \eqref{estimbarch} with $\bar{C}_c = C_{tr} \tilde{C}(\eta,\tilde{d})$.     
\QED \\

As a consequence of Propositions \ref{coercivebarbh} and \ref{boundbarch} we have
\begin{theorem}
\label{coercivebarah}
Provided $h$ is small enough, the bilinear form $\bar{a}_h$ is uniformly coercive over both $\tilde{V}^k_h$ (resp. $V_h^k$), in the sense that there exists a mesh-independent constant $\bar{\alpha}^{'} > 0$ such that 
\begin{equation}
\label{coercivbarah}
\bar{a}_h(v,v) \geq \bar{\alpha}^{'} \| v \|_{1,h}^2 \; \forall v \in \tilde{V}^k_h \mbox{ (resp. } \forall v \in V_h^k\mbox{).} 
\end{equation}
\end{theorem}
\prov First we observe that \eqref{estimbarch} trivially extends to $v \in \tilde{V}^k_h$. Thus taking $w=v \in \tilde{V}_h^k$ we have 
$\bar{c}_h(v,v) \geq - \bar{C}_c h^{1/2} \| v \|_{1,h}^2$. Taking into account \eqref{coercivbarbh} it is readily seen that, provided 
$h \leq [\bar{\beta}/(2 \bar{C}_c)]^2$ \eqref{coercivbarah} holds with $\bar{\alpha}^{'} = \bar{\beta}/2$. \QED \\

Before switching to the bilinear form $a_h$ we prove some auxiliary results:
\begin{lemma}
\label{auxiliary1}
Let $R_h^k$ be the bilinear form defined on $V^k_h \times V^k_h$ by   
\begin{equation}
\label{Rhk}
\left\{
\begin{array}{l}
R_h^k(v,w) = \displaystyle \sum_{T \in {\mathcal T_h}} R_T(w,v) \; \forall (v;w) \in [V^k_h]^2, \\
\mbox{where}\\
R_T^k(v,w) = \int_{T} \tilde{r} w v - J_T^k[\tilde{r}wv].
\end{array}
\right.
\end{equation}
If the quadrature formula \eqref{JTk} is exact for polynomials of degree less than or equal to two  
and $r_0 > 0$, there exists a mesh-independent constant $C(r) > 0$ such that it holds
\begin{equation}
\label{lemmaux1} 
R_h^k(w,v) \leq C(r) h \| w \|_{1,h} \| v \|_{1,h} \; \forall (v;w) \in [V^k_h]^2. 
\end{equation}
\end{lemma}
\prov
Extending eqn. (2.39) in Theorem 2.1 of \cite{Nedoma} to functions $\varphi \in H^m(T)$ $\forall T \in {\mathcal T}_h$, setting $\phi=(\tilde{r} w v)_{|T}$ and $m=2$ we obtain for a constant $C_{\Omega}$ independent of $h_T$
\begin{equation}
\label{lemmaux11}
R_T^k(\tilde{r} w v) \leq C_{\Omega} h_T^2 \| \tilde{r} w \|_{2,T} \| v \|_{1,T} \; \forall T \in {\mathcal T}_h.
\end{equation}
It turns out that, for a suitable integer $J$, we have $\forall T \in {\mathcal T}_h$ 
\begin{equation}
\label{lemmaux12}
\| \tilde{r} w \|_{2,T} \leq J [\| \tilde{r} \|_{0,\infty,T} \| w \|_{2,T} + \| \nabla \tilde{r} \|_{0,\infty,T} \| w \|_{1,T} + \| H(\tilde{r}) \|_{0,\infty,T} \| w \|_{0,T}], 
\end{equation}
assuming that $\tilde{r} \in \displaystyle C^1(\overline{\tilde{\Omega}_h})$. \\
On the other hand, using a classical inverse inequality (cf. \cite{Verfuerth}), we know that $\| w \|_{2,T} \leq C_{\iota} h_T^{-1} 
\| w \|_{1,T}$ for a constant $C_{\iota}$ independent of $T$. Plugging this inequality into \eqref{lemmaux12} and then the resulting inequality into \eqref{lemmaux11}, from the above definition of $R_h$ we immediately obtain \eqref{lemmaux1} for a constant $C(r)$ expressed in terms of $\| \tilde{r} \|_{2,\infty}$ . \QED   

\begin{lemma}
\label{auxiliary2}
Let $D_h^k$ be the bilinear form defined on $V^k_h \times V^k_h$ by   
\begin{equation}
\label{Dhk}
\left\{
\begin{array}{l}
D_h^k(v,w) = \displaystyle \sum_{T \in {\mathcal T_h}} D^k_T(w,v) \; \forall (v;w) \in [V^k_h]^2, \\
\mbox{where}\\
D_T^k(v,w) = \int_{T} \tilde{d} \nabla w \cdot \nabla v - J_T^k[\tilde{d} \nabla w \cdot \nabla v].
\end{array}
\right.
\end{equation}
Provided the numerical quadrature formula $J^k_T(\phi)$ is exact whenever $\phi \in P_{2k-2}(T)$, there exists a mesh-independent constant $C(d) > 0$ such that it holds

\begin{equation}
\label{lemmaux2} 
D_h^k(w,v) \leq C(d) h \| w \|_{1,h} \| v \|_{1,h} \; \forall (v;w) \in [V^k_h]^2. 
\end{equation}
\end{lemma}
\prov
\eqref{lemmaux2} is a direct application of Theorem 2 in \cite{Zenisek} in the particular case where the diffusion tensor equals $d$ multiplied by the identity tensor and the mapping from any current triangle $T \in {\mathcal T}_h$ onto the master triangle $\hat{T}$ is affine. Such a theorem in turn is based on references given in \cite{Zenisek} itself and corroborated by other results available in the literature such as those in \cite{BarrettElliott}. \QED  \\

Now let us denote by ${\bf t}_h$ the unit tangent vector to $\Gamma_h$ defined by a rotation of ${\bf n}_h$ by $\pi/2$ in the counterclockwise sense. We represent by $\mbox{d}\varphi/\mbox{d}\tau_h$ the first order derivative along $\Gamma_h$ in the direction of ${\bf t}_h$ of a differentiable function $\varphi$ defined on $e_T$ for an arbitrary $T \in {\mathcal S}_h$. \\
The following result is a more general variant of Lemma 3.31 of \cite{FeistauerNajzar} in the case of hilbertian norms.
\begin{lemma}
\label{auxiliary4}
For a given $T \in {\mathcal S}_h$ and $j \geq 0$, we define a bounded linear form ${\mathcal H}_T^j$ on $W^{j+3/2,1}(T)$ by ${\mathcal H}_T^j(\cdot) := \int_{e_T} (\cdot) - {\mathcal I}^k_T(\cdot)$. Let $\phi$ be a given function in $C^{j+1}(e_T)$. Provided the integration formula ${\mathcal I}^k_T$ is exact for polynomials of degree less than or equal to $j$ in $e_T$, there exists a constant $C_{\mathcal H}$ independent of $T$ such that $\forall z \in H^{j+3/2}(T)$ and  $\forall v  \in P_k(T)$ it holds for $n:=\min[k,j+1]$
\begin{equation}
\label{lemmaux41}
\left\{
\begin{array}{l}
{\mathcal H}_T^j(\phi z v) \leq C_{\mathcal H} \displaystyle \left[ h_T^{j+1} \| v \|_{0,e_T} \displaystyle \left(\| z \|_{j+1,e_T} \| \phi \|_{0,\infty,e_T} + \displaystyle \sum_{l=1}^{j+1} | z |_{j+1-l,e_T} | \phi |_{l,\infty,e_T} \right) \right. \\
\left. + h_T^{j-n+3/2} \| v \|_{1/2,e_T} \displaystyle \sum_{l=0}^{j} \| z \|_{j-l,e_T}) \| \phi \|_{l,\infty,e_T} \right]. 
\end{array}
\right.
\end{equation}
\end{lemma} 
\prov
First we note that ${\mathcal H}_T^j$ is a bounded linear functional on $W^{j+1,1}(e_T)$ for any $j \geq 0$. Thus, resorting to the master element $\hat{T}$ and defining   
\begin{equation}
\label{lemmaux42}
\hat{\mathcal H}^j(\hat{\varphi}) = \int_{\hat{e}} \hat{\varphi} - \hat{\mathcal I}^k(\hat{\varphi}) \; \forall \hat{\varphi} \in W^{j+1,1}(\hat{e}),
\end{equation}
we observe that
\begin{equation}
\label{lemmaux43}
{\mathcal H}_T^j(\varphi) = |e_T| \hat{{\mathcal H}}^j(\hat{\varphi})/\hat{e} \; \forall \varphi \in W^{j+1,1}(e_T).
\end{equation} 
Clearly enough, $\hat{{\mathcal H}}^j$ is a uniformly bounded linear functional on $W^{j+1,1}(\hat{e})$ for any $j \geq 0$, independently of $T$. 
Furthermore, from the assumptions on ${\mathcal H}_T^j$, we can assert that $\hat{\mathcal H}^j(\hat{\varphi})=0$ whenever $\hat{\varphi} \in P_{j+1}(\hat{e})$. Therefore, from the Bramble-Hilbert Lemma, there exists a constant $\hat{C}_{\mathcal H}$ such that
$$\hat{{\mathcal H}}^j(\hat{\varphi}) \leq \hat{C}_{\mathcal H} | \hat{\varphi} |_{j+1,1,\hat{e}} \; \forall \hat{\varphi} \in H^{j+1}(\hat{e}).$$
Plugging the above upper bound into \eqref{lemmaux43} and moving back to the current triangle $T \in {\mathcal S}_h$, from classical results we obtain for a suitable constant $C^{'}_{\mathcal H}$ independent of $T$ 
\begin{equation}
\label{lemmaux44}
{\mathcal H}_T^j(\varphi) = C^{'}_{\mathcal H} h_T^{j+1} | \varphi |_{j+1,1,e_T} \; \forall \varphi \in W^{j+1,1}(e_T).
\end{equation}   
Now setting $\varphi = \phi z v$ in \eqref{lemmaux44}, by straightforward calculations, we successively come up with two other constants $C^{''}_{\mathcal H}$ and $C^{'''}_{\mathcal H}$ independent of $T$ such that  $ \forall z \in H^{j+3/2}(T)$ and $\forall v \in P_k(T)$ it holds
\begin{equation}
\label{lemmaux45}
{\mathcal H}_T^j(\phi zv) = C^{'''}_{\mathcal H} h_T^{j+1} (\| \phi \|_{0,\infty,e_T} | zv |_{j+1,1,e_T} + \displaystyle \sum_{l=0}^{j}\|\mbox{d}\phi/\mbox{d}\tau \|_{l,\infty,e_T} \| zv \|_{j-l,1,e_T}. 
\end{equation}
Noticing that $| v |_{m,e_T} = 0$ for $m > k$, it also holds for another mesh-independent constant $\bar{C}_{\mathcal H}$ 
\begin{equation}
\label{lemmaux46}
\left\{
\begin{array}{l}
{\mathcal H}_T^j(\phi zv) = \bar{C}_{\mathcal H} h_T^{j+1} \displaystyle \left[\| v \|_{0,e_T} \displaystyle \left(\| z \|_{j+1,e_T} \| \phi \|_{0,\infty,e_T}  + \displaystyle \sum_{l=1}^{j+1} | z |_{j+1-l,e_T} | \phi |_{l,\infty,e_T})\right) \right. \\
\left. + \displaystyle \left\| \frac{\mbox{d}v}{\mbox{d}\tau_h} \right\|_{n-1,e_T} \displaystyle \sum_{l=0}^{j} \| z \|_{j-l,e_T}) \| \phi \|_{l,\infty,e_T} \right]
\end{array}
\right.
\end{equation} 
Since there exists a constant $C_{n,1/2}$ independent of $T$ such that  $\| \mbox{d}v/\mbox{d}\tau_h \|_{n-1,e_T} \leq C_{n,1/2} h_T^{1/2-n} \| v \|_{1/2,e_T}$ (cf. \cite{Georgoulis}), the result follows. \QED
     
\begin{e-proposition}
\label{coercivebh}
Let the integration formula \eqref{JTk} be exact for polynomials of degree at least equal to $2k-2$. Provided $h$ is sufficiently small, the bilinear form $b_h$ is uniformly coercive over $V^k_h$ for the $\| \cdot \|_{1,h}$-norm, in the sense that there exists a constant $\beta > 0$ independent of $h$ such that 
\begin{equation}
\label{beta}
b_h(v,v) \geq \beta \| v \|_{1,h}^2 \; \forall v \in V_h^k.
\end{equation}
\end{e-proposition}  

\prov
First of all we apply Lemma \ref{auxiliary4} to the term of $b_h(w,v)$ equal to the summation of ${\mathcal I}_T^k(\phi w v)$ over ${\mathcal S}_h$ for $\phi = \mbox{cos } \theta_T \bar{\eta}$, with $j=0$. Since the sum of the weights $\omega_i$ equals one, we may write
\begin{equation}
\label{coercivbh0}
{\mathcal H}^0_T( \phi w v) \leq C_{\mathcal H} [ h_T \| v \|_{0,e_T} ( \| \phi \|_{0,\infty,e_T} | w |_{1,e_T}  + | \phi |_{1,\infty,e_T} \| w \|_{0,e_T}) + h_T^{1/2} \| v \|_{1/2,e_T} \| \phi \|_{0,\infty,e_T} \| w \|_{0,e_T}].
\end{equation}
Now we trivially have 
\begin{equation}
\label{coercivbh1}
\| \phi \|_{0,\infty,e_T} \leq \| \eta \|_{0,\infty,\Gamma} \; \forall T \in {\mathcal S}_h. 
\end{equation}
On the other hand, referring to \cite{ZAMM}, we know that $|\mbox{d}(\mbox{cos }\theta_T)/\mbox{d}\tau_h| = |\mbox{sin }\theta_T \mbox{d}\theta_T/\mbox{d}\tau_h| \leq C^{'}(\Gamma) s_m h_T $ where $C^{'}(\Gamma)$ depends only on $\Gamma$. Furthermore, again according to \cite{ZAMM}, $|\mbox{d}\bar{\eta}/\mbox{d}\tau_h| \leq (1+C_{\kappa} h_T)\mbox{d}\eta/\mbox{d}\tau$, where $\mbox{d}(\cdot)/\mbox{d}\tau$ represents the tangential derivative along $\Gamma$ in the direction of ${\bf t}$ and $C_{\kappa}$ is a constant independent of $T$. It follows that
\begin{equation}
\label{coercivbh2}
| \phi |_{1 ,\infty,e_T} \leq C^{'}(\Gamma) s_m h_T \| \eta \|_{0 ,\infty,\Gamma} + (1+C_{\kappa} h_T) |\eta|_{1,\infty,\Gamma}, 
\end{equation}
or yet, for a a mesh-independent constant $C^{''}_{\Gamma}$,
\begin{equation}
\label{coercivbh2bis}
| \phi |_{1 ,\infty,e_T} \leq C^{''}(\Gamma) \| \eta \|_{1 ,\infty,\Gamma}. 
\end{equation} 
Plugging \eqref{coercivbh1} and \eqref{coercivbh2bis} into \eqref{coercivbh0}
, we readily obtain
\begin{equation}
\label{coercivbh4}  
{\mathcal H}^0_T( \mbox{cos } \theta_T \bar{\eta} w v) \leq C^{''}(\eta) [h_T \| w \|_{1,e_T} \| v \|_{0,e_T} + h_T^{1/2} \| w \|_{0,e_T} \| v \|_{1/2,e_T}],
\end{equation}
where 
\begin{equation}
\label{C2primeeta}
C^{''}(\eta) = C_{\mathcal H} \sqrt{(C_{\Gamma}^{''})^2+1} \| \eta \|_{1,\infty,\Gamma}.
\end{equation}
Using the inverse inequality $\| w \|_{1,e_T} \leq C_{1,1/2} h_T^{-1/2} \| w \|_{1/2,e_T}$ for $w \in P_k(T)$, \eqref{coercivbh4} can be shortened into 
\begin{equation}
\label{coercivbh5}  
{\mathcal H}^0_T( \mbox{cos } \theta_T \bar{\eta} w v) \leq C^{'}(\eta) h_T^{1/2} \| w \|_{1/2,e_T} \| v \|_{1/2,e_T},
\end{equation}
for another constant $C^{'}(\eta)$ depending only on $\Gamma$ and $\eta$.\\ 
Taking $w=v$, summing up over ${\mathcal S}_h$ and using the trace inequality \eqref{trace} we readily obtain   
\begin{equation}
\label{coercivbh6}  
\displaystyle \sum_{T \in {\mathcal S}_h} {\mathcal H}^0_T( \mbox{cos } \theta_T \bar{\eta} v^2) \leq C(\eta) h_T^{1/2} \| v \|^2_{1,T},
\end{equation}
with $C(\eta) = C^{'}(\eta) C_{tr}^2$.\\
\indent Now let us define $C_0(\eta) = C(\eta)$ if $\eta_0 >0$ and $C_0(\eta)=0$ otherwise, together with $C_0(r) = C(r)$ if $r_0 >0$ and $C_0(r)=0$ otherwise. \\  
Given $v \in \tilde{V}_h^k$ we apply 
Lemmata \eqref{auxiliary1} and \eqref{auxiliary2} with $w = v$ together \eqref{coercivbh6}, to readily obtain  
\begin{equation}
\label{coercivbh1bis}
\left\{
\begin{array}{l}
b_h(v,v) \geq \bar{b}_h(v,v) - {\bf C}(h) \| v \|_{1,h}^2 \\
\mbox{where}\\
{\bf C}(h) := C(d)h + C_0(\eta) h^{1/2} + C_0(r) h.
\end{array}
\right.
\end{equation} 
Thus, taking into account \eqref{coercivbarbh}, as long as $h$ is such that ${\bf C}(h) \leq \bar{\beta}/2$, \eqref{beta} holds with 
$\beta=\bar{\beta}/2$. \QED \\

We proceed by addressing the bilinear form $c_h$ as follows.\\
\begin{e-proposition}
\label{coercivech}
There exists a mesh-independent constant $C_c$ such that
\begin{equation}
\label{lowerboundch}
c_h(v,v) \geq - C_c h^{1/2} \| v \|_{1,h}^2.
\end{equation}
\end{e-proposition}

\prov
Akin to the proof of Proposition \ref{coercivebh}, we resort to Lemma \ref{auxiliary4} combined with some additional calculations. However here we must handle three terms instead of just one. More precisely, referring to the definition of the functional ${\mathcal H}_T^j$ in the statement of Lemma \ref{auxiliary4}, this means that we have to estimate the difference between $\bar{c}_h(v,v) - c_h(v,v)$ for $v \in \tilde{V}^k_h$, by properly bounding above the summation over ${\mathcal S}_h$ of the following terms:
\begin{itemize}
\item $\delta_T^1(w,v):= {\mathcal H}_T^0(\mbox{cos } \theta_T \bar{\eta} [\bar{w} - w]v)$;
\item $\delta_T^2(w,v):=\displaystyle {\mathcal H}_T^0\left(\mbox{cos }\theta_T \left[ \displaystyle \int_M^P \frac{\partial \left(\tilde{d} \displaystyle \frac{\partial w}{\partial n_P}\right)}{\partial n_h} \right] v \right)$;
\item $\delta_T^3(w,v):= \displaystyle {\mathcal H}_T^0\left(\mbox{sin } \theta_T  \tilde{d} \frac{\partial w}{\partial \tau_P} \; v \right)$.
\end{itemize}
In the case of $\delta_T^1$ we may apply Lemma \ref{auxiliary4} with $\phi = \mbox{cos } \theta_T \bar{\eta}$ and $z = \int_M^P \partial w/\partial n_h$. In doing so, recalling \eqref{C2primeeta}, it holds,
\begin{equation}
\label{coercivch11}
\delta_T^1(w,v)\leq C^{''}(\eta) \displaystyle \left( h_T  \left| \displaystyle \int_M^P \frac{\partial w}{\partial n_h} \right|_{1,e_T}  \| v \|_{0,e_T} +  h_T^{1/2} \left\| \int_M^P \frac{\partial w}{\partial n_h} \right\|_{0,e_T}  \| v \|_{1/2,e_T} \right).
\end{equation}
Noticing that $h_T \leq 1$, we assemble both terms on the right hand side of  \eqref{coercivch11} to come up with another constant $C^{''}_1(\eta)$ such that
\begin{equation}
\label{coercivch12}
\delta_T^1(w,v)\leq C^{''}_1(\eta) h_T^{1/2} \displaystyle \left\|  \displaystyle \int_M^P \frac{\partial w}{\partial n_h} \right\|_{1,e_T}  \| v \|_{1/2,e_T}. 
\end{equation}
In order to properly estimate the third factor on the right hand side of \eqref{coercivch12} we make use of the following identity applying to any twice differentiable function $\varphi$ in $\tilde{T}$.  
\begin{equation}
\label{dtauh}
\left\{
\begin{array}{l}
\displaystyle \frac{ \mbox{d} \displaystyle \int_M^P \left( \frac{\partial \varphi}{\partial n_h}\right)}{\mbox{d}\tau_h} = \displaystyle \left[\frac{\mbox{d}\varphi}{\mbox{d}\tau_h}\right]\!(P) - \displaystyle \left[\frac{\partial \varphi}{\partial \tau_h}\right]\!(M) \\
= \displaystyle \left[\frac{\partial \varphi}{\partial \tau_h}\right]\!(P) + \left[\frac{\partial \varphi}{\partial n_h}\right]\!(P) \displaystyle \frac{\mbox{d}\overline{MP}}{\mbox{d}\tau_h} - \displaystyle \left[\frac{\partial \varphi}{\partial \tau_h}\right]\!(M) \\ 
= \displaystyle \left[\frac{\partial \varphi}{\partial n_h}\right]\!(P)  \mbox{ tan }\theta_T 
+ \displaystyle \int_M^P \frac{\partial^2 \varphi}{\partial n_h \partial \tau_h}. 
\end{array}
\right.
\end{equation} 
Taking $\varphi = w$ in \eqref{dtauh} it is readily seen that
\begin{equation} 
\label{coercivch12bis}
\left\{
\begin{array}{l}
\displaystyle \left| \int_M^P \frac{\partial w}{\partial n_h} \right|_{1,e_T} \leq  \displaystyle \frac{s_m}{c_m} h_T^{3/2} \displaystyle \left\| \frac{\partial w}{\partial n_h} \right\|_{0,\infty,\tilde{T}} + C_{\Gamma}^{1/2} h_T | w |_{2,\Delta_T}\\
\mbox{whereas}\\
\| \int_M^P \int_M^P \partial w/\partial n_h \|_{0,e_T} \leq C_{\Gamma}^{1/2} h_T | w |_{1,\Delta_T}.
\end{array}
\right.
\end{equation}
Since there exists a mesh-independent constant $\tilde{C}_I$ such that (cf. \cite{ZAMM}) 
$$\| \partial w/\partial n_h \|_{0,\infty,\tilde{T}} \leq \tilde{C}_I h_T^{-1} | w |_{1,T}, $$
taking also into account a standard inverse inequality relating $\| \cdot \|_{2,T}$ to $\| \cdot \|_{1,T}$, 
after plugging \eqref{coercivch12bis} into \eqref{coercivch12} we come up with a constant 
$C^{'}_1(\eta)$ depending only on $\Gamma$ and $\eta$ such that  
\begin{equation} 
\label{coercivch13}
\delta_T^1(w,v)  \leq C^{'}_1(\eta) h_T^{1/2} \| w \|_{1,T} \| v \|_{1/2,e_T}. 
\end{equation}
Sweeping ${\mathcal S}_h$, using the trace inequality \eqref{trace} and setting $C_1(\eta) = C^{'}_1(\eta) C_{tr}$, the following upper bound immediately derives from \eqref{coercivch13}.
\begin{equation}
\label{coercivch14}
\displaystyle \sum_{T \in {\mathcal S}_h} \delta_T^1(w,v) \leq C_1(\eta) h^{1/2} \| w \|_{1,h} \| v \|_{1,h}. 
\end{equation}
In order to better fit the bilinear form $\delta_T^2$ into the functional framework of Lemma \ref{auxiliary4}, it is judicious to split it 
into the sum of two bilinear forms $\delta_T^{2a}$ and $\delta_T^{2b}$ given by  
\begin{equation}
\label{coercivch2a}
\delta_T^{2a}(w,v) := \displaystyle {\mathcal H}_T^0\left(\mbox{cos}^2 \theta_T \displaystyle \left[  \displaystyle \int_M^P \frac{\partial \left(\tilde{d} \displaystyle \frac{\partial w}{\partial n_h}\right)}{\partial n_h} \right] v \right)
\end{equation} 
and
\begin{equation}
\label{coercivch2b}
\delta_T^{2b}(w,v) := -\displaystyle {\mathcal H}_T^0\left(\mbox{cos } \theta_T \mbox{sin } \theta_T \displaystyle \left[  \displaystyle \int_M^P \frac{\partial \left(\tilde{d} \displaystyle \frac{\partial w}{\partial \tau_h}\right)}{\partial n_h} \right] v \right).
\end{equation}
In the case of $\delta_T^{2a}$ we have $\phi = \mbox{cos}^2 \theta_T$ and $z = \displaystyle \int_M^P \displaystyle \frac{\partial( \tilde{d} \partial w/\partial n_h)}{\partial n_h}$. \\
Using elementary Calculus,    
on the one hand we have $$\| \phi \|_{1,\infty,e_T} \leq 1+2 s_m h_T \| \mbox{d}\theta_T/\mbox{d}\tau_h \|_{0,\infty,e_T},$$ i.e., $\exists C^{'}_{2a}(\Gamma)$ independent of $T$ such that 
\begin{equation}
\label{boundphi2a}
\| \phi \|_{1,\infty,e_T} \leq C^{'}_{2a}(\Gamma) \mbox{ with } \phi = \mbox{cos}^2 \theta_T.
\end{equation}
On the other hand, taking $\varphi = \tilde{d} \partial w/\partial n_h$ in \eqref{dtauh} and $z = \int_M^P \partial \varphi/\partial n_h$, after straightforward calculations we obtain    
\begin{equation} 
\label{coercivch21}
\left\{
\begin{array}{l}
\| z \|_{0,e_T} \leq \| \tilde{d} \|_{1,\infty,\tilde{\Omega}_h} C_{\Gamma}^{1/2} h_T \displaystyle \|w\|_{2,\Delta_T}  \\
\mbox{and} \\
| z |_{1,e_T} \leq \| \tilde{d} \|_{1,\infty,\tilde{\Omega}_h} h_T \left[ \displaystyle \frac{s_m}{c_m} \displaystyle \left( \left\| \frac{\partial w}{\partial n_h} \right\|_{0,e_T} + \displaystyle \left\| \frac{\partial^2 w}{\partial n_h^2} \right\|_{0,e_T}\right) \right. \\
\left. + C_{\Gamma}^{1/2} \left( | w |_{2,\Delta_T} + | w |_{3,\Delta_T} \right) \right].
\end{array}
\right.  
\end{equation}
Using inequalities such as \eqref{GStrangGFix} already repeatedly exploited in this work, we can successively write for two pairs of mesh-independent constants $(C^{'}_0(d);C^{'}_1(d))$ and $(C^{''}_0(d);C^{''}_1(d))$
\begin{equation} 
\label{coercivch22}
\left\{
\begin{array}{l}
\| z \|_{0,e_T} \leq C^{'}_0(d) h_T^{3/2} \|w\|_{2,T} \leq C^{''}_0(d) h_T^{1/2} \|w\|_{1,T} \\
\mbox{and} \\
| z |_{1,e_T} \leq  C^{'}_1(d) h_T^{3/2} (\|w\|_{2,\infty,T} + \| w \|_{3,T}) \leq C^{''}_1(d) h_T^{-1/2} \|w\|_{1,T}.
\end{array}
\right.  
\end{equation} 
Now we apply Lemma \ref{auxiliary4} with $j=0$, by plugging \eqref{coercivch22} and \eqref{boundphi2a} into \eqref{lemmaux41}. This leads to a constant $C_{2a}^{'}(d)$ depending only on $d$ and $\Gamma$ such that   
\begin{equation} 
\label{coercivch23a}
\delta_T^{2a}(w,v) \leq C_{2a}^{'}(d) \| w \|_{1,T} (h_T^{1/2} \| v \|_{0,e_T} + h_T  \| v \|_{1/2,e_T}).
\end{equation}

In the case of $\delta_T^{2b}(w,v)$ we have again $z = \int_M^P \partial \varphi/\partial n_h$ though with $\varphi = \tilde{d} \partial w/\partial \tau_h$. Moreover, since now $\phi = \mbox{cos }\theta_T \mbox{sin }\theta_T$, the estimate of this term using \eqref{dtauh} is nothing but a more favorable variant of the previous one. Indeed, with no difficulty, we determine a constant $C_{2b}^{'}(d)$ depending only on $d$ and $\Gamma$ such that   
\begin{equation} 
\label{coercivch23b}
\delta_T^{2b}(w,v) \leq C_{2b}^{'}(d) \| w \|_{1,T} (h_T^{1/2} \| v \|_{0,e_T} + h_T  \| v \|_{1/2,e_T}).
\end{equation}
Combining \eqref{coercivch23a} and \eqref{coercivch23b}, similarly to the case of $\delta_T^{1}(w,v)$, we obtain for a constant 
$C_2(d)$ depending only on $\Gamma$ and $d$
\begin{equation}
\label{coercivch24}
\displaystyle \sum_{T \in {\mathcal S}_h} \delta_T^2(w,v) \leq C_2(d) h^{1/2} \| w \|_{1,h} \| v \|_{1,h}. 
\end{equation}
Finally we address the case of $\delta_T^3(w,v)$, which can also by split into the sum of $\delta_T^{3a}(w,v)$ and $\delta_T^{3b}(w,v)$ given by
\begin{equation}
\label{coercivch3a}
\delta_T^{3a}(w,v) := \displaystyle {\mathcal H}_T^0\left(\mbox{sin } \theta_T \mbox{cos } \theta_T \tilde{d} \displaystyle \left[ \frac{\partial w}{\partial \tau_h} \right] v \right)
\end{equation} 
and
\begin{equation}
\label{coercivch3b}
\delta_T^{3b}(w,v) := \displaystyle {\mathcal H}_T^0\left(-\mbox{sin}^2 \theta_T \tilde{d} \displaystyle \left[ \frac{\partial w}{\partial n_h} \right] v \right)
\end{equation}
In the case of $\delta_T^{3a}(w,v)$ we take $\phi=\mbox{sin } \theta_T \mbox{cos } \theta_T \tilde{d}$ and $z=\partial w/\partial \tau_h$. It is easily seen that there exist two constants $C {'}_{3a}(=s_m)$ and $C^{''}_{3a}$ depending only on $\Gamma$ such that 
\begin{equation}
\label{boundphi3a}
\left\{
\begin{array}{l}
\| \phi \|_{0,\infty,e_T} \leq  C^{'}_{3a} h_T \| \tilde{d} \|_{0,\infty,\tilde{\Omega}_h} \\
\mbox{and}\\
| \phi |_{1,\infty,e_T} \leq  C^{''}_{3a} \| \tilde{d} \|_{1,\infty,\tilde{\Omega}_h}.
\end{array}
\right.
\end{equation}
Moreover, we have
\begin{equation} 
\label{coercivch31}
\left\{
\begin{array}{l}
\| z \|_{0,e_T} \leq | w |_{1,e_T} \\
\mbox{and}\\
| z |_{1,e_T} \leq | w |_{2,e_T}. 
\end{array}
\right.  
\end{equation} 
Using \eqref{boundphi3a} and \eqref{coercivch31}, application of Lemma \ref{auxiliary4} with $j=0$, assorted with some inverse inequalities already repeatedly exploited in this work, directly yields for a suitable mesh-independent constant $C_{3a}(\tilde{d})$ 
\begin{equation} 
\label{coercivch3a1}
\delta_T^{3a}(w,v) \leq C_{3a}(\tilde{d}) h_T^{1/2} (\| w \|_{1/2,e_T} \| v \|_{0,e_T} +  \| w \|_{0,e_T} \| v \|_{1/2,e_T}) 
\end{equation}
As for $\delta_T^{3b}(w,v)$, we have $\phi = \mbox{ sin}^2 \theta_T$ and $z = \partial w/\partial n_h$. Therefore, quite obviously, for two constants $C_{3b}^{'}$ and $C_{3b}^{''}$ independent of ${\mathcal T}_h$ it holds
\begin{equation}
\label{boundphi3b}
\left\{
\begin{array}{l}
\| \phi \|_{0,\infty,e_T} \leq  C^{'}_{3b} h_T^2 \| \tilde{d} \|_{0,\infty,\tilde{\Omega}_h}, \\
\mbox{and}\\
| \phi |_{1,\infty,e_T} \leq  C^{''}_{3b} h_T \| \tilde{d} \|_{1,\infty,\tilde{\Omega}_h}
\end{array}
\right.
\end{equation}
and furthermore
\begin{equation} 
\label{coercivch33}
\left\{
\begin{array}{l}
\| z \|_{0,e_T} \leq h_T^{1/2} \| \nabla w \|_{0,\infty,T} \\
\mbox{and}\\
| z |_{1,e_T} \leq h_T^{1/2} \| H(w) \|_{0,\infty,T},  
\end{array}
\right.  
\end{equation}
which yields for another constant $C^{'''}_{3b}$ independent of $T$
\begin{equation} 
\label{coercivch34}
\left\{
\begin{array}{l}
\| z \|_{0,e_T} \leq h_T^{-1/2} \| w \|_{1,T} \\
\mbox{and}\\
| z |_{1,e_T} \leq h_T^{-3/2} \| w \|_{1,T}  
\end{array}
\right.  
\end{equation}
Applying Lemma \eqref{auxiliary4} with $j=0$, from \eqref{boundphi3b} and \eqref{coercivch34} we infer the existence of a mesh-independent constant  
$C_{3b}(\tilde{d})$ such that
\begin{equation} 
\label{coercivch3b1}
\delta_T^{3b}(w,v) \leq C_{3b}(\tilde{d}) | w |_{1,T} ( h_T^{3/2}  \| v \|_{0,e_T} +  h_T \| v \|_{1/2,e_T}) 
\end{equation}
Therefore, taking into account \eqref{coercivch3a1} and summing up over ${\mathcal S}_h$, by the trace inequality \eqref{trace} it trivially follows that that there is a constant $C_3(\tilde{d})$ independent of ${\mathcal T}_h$ such that
\begin{equation}
\label{coercivch35}
\displaystyle \sum_{T \in {\mathcal S}_h} \delta_T^3(w,v) \leq C_3(\tilde{d}) h^{1/2} \displaystyle \| w \|_{1,h} \| v \|_{1,h} 
\end{equation} 
Finally, combining \eqref{coercivch14}, \eqref{coercivch24} and \eqref{coercivch35}, we establish
\begin{equation}
\label{chwv6} 
- c_h(w,v) \leq C_c h^{1/2} \| w \|_{1,h} \| v \|_{1,h} \; \forall w \in \tilde{V}^k_h \mbox{ and } \forall v \in V^k_h,
\end{equation}
with $C_C = C_1(\eta)+C_2(\tilde{d})+C_3(\tilde{d})$. Taking $w=v$ in $\Omega_h$, this yields \eqref{lowerboundch}.     
\QED \\

Now, the following theorem immediately derives from Propositions \ref{coercivebh} and \ref{coercivech}.
\begin{theorem}
\label{coerciveah}
There exists a constant $\alpha^{'} > 0$ independent of $h$ such that for every $v \in \tilde{V}^k_h$ (resp. $v \in V_h^k$) it holds 
\begin{equation}
\label{coercivah}
a_h(v,v) \geq \alpha^{'} \| w \|_{1,h}^2   
\end{equation}
\end{theorem}
\prov
Taking into account \eqref{ah}, \eqref{bh}, \eqref{ch} together with \eqref{beta} and \eqref{lowerboundch}, as long as $h$ is small enough, and in any case $h^{1/2} \leq \beta / (2 C_c)$, \eqref{coerciveah} holds with 
$\alpha{'} = \beta/2$. \QED \\

\subsection{Uniform weak-coercivity}
 
Theorem \ref{coerciveah} establishes the uniform coercivity of the bilinear form $a_h$ over $\tilde{V}_h^k$ (resp. $V_h^k$) equipped with the norm $\| \cdot \|_{1,h}$. However, in the sequel we will need a variant of this result, namely the uniform weak coercivity of $a_h$ over $\tilde{V}_h^k \times V_h^k$, where $\tilde{V}_h^k$ is equipped with an extended norm defined as below.\\
Recalling the set $\tilde{\Omega}_h$ defined in Subsection 3.1, we use the special notation $\| \cdot \|_{j,\tilde{h}}$ for the norm of $H^j(\tilde{\Omega}_h)$, and extend it to the broken $H^j$-norm for a given $j>1$, i.e., the natural norm of the space $\tilde{H}^j(\tilde{\Omega}_h)$ consisting of functions $\tilde{w} \in H^1(\tilde{\Omega}_h)$ whose restriction to $\tilde{T}$ for every $T \in {\mathcal T}_h$ belongs to $H^j(\tilde{T})$. This means that 
$$ \| \tilde{w} \|_{j,\tilde{h}} := \displaystyle \left[ \sum_{T \in {\mathcal T}_h} \| \tilde{w} \|_{j,\tilde{T}}^2 \right]^{1/2} \; \forall \tilde{w} \in \tilde{H}^j(\tilde{\Omega}_h).$$
Let $\tilde{W}^j_h$ be the direct sum of $H^j(\tilde{\Omega}_h)$ and $\tilde{V}^k_h$, actually a subspace of $\tilde{H}^j(\tilde{\Omega}_h)$.
For the space $\tilde{W}^2_h$ we introduce the following mesh-dependent norm: 
\begin{equation}
\label{NormtildeW2h}
\| \tilde{w} \|_{\tilde{W},2,h} := \| \tilde{w} \|_{1,h} + h \| \tilde{w} \|_{2,\tilde{h}} \; \forall \tilde{w} \in \tilde{W}^2_h.
\end{equation}
It is not difficult to figure out that this norm is uniformly equivalent to $\| \cdot \|_{1,h}$ over $\tilde{V}^k_h$, in the sense that there exist two mesh-independent constants $C^{-}$ and $C^{+}$ such that 
\begin{equation}
\label{equiv}
C^{-} \| w \|_{\tilde{W},2,h} \leq \| w \|_{1,h} \leq C^{+} \| w \|_{\tilde{W},2,h} \; \forall w \in \tilde{V}^k_h, 
\end{equation}  
Indeed, while the right inequality trivially holds with $C^{+} =1$, the left inequality directly follows from standard inverse inequalities for the norms of $H^j(T)$ or $H^j(\tilde{T})$ (cf. \cite{ErnGuermond,ZAMM}).\\
From this observation follows    
\begin{theorem} 
\label{weakcoercivity}
Both bilinear forms $\bar{a}_h$ and $a_h$ are uniformly weakly coercive on the product of the spaces $\tilde{V}_h^k$ and $V_h^k$ respectively equipped with the norms $\| w \|_{\tilde{W},2,h}$ and $\| \cdot \|_{1,h}$, that is, there exists mesh-independent constant $\bar{\alpha} >0$ and $\alpha >0$ such that
\begin{equation}
\label{infsup}
\left\{
\begin{array}{l}
\forall w \in \tilde{V}_h^k \; \displaystyle \sup_{v \in V_h^h \setminus \{0\}} \frac{\bar{a}_h(w,v)}{\| v \|_{1,h}} \geq \bar{\alpha} \| w \|_{\tilde{W},2,h} \\
\mbox{and}\\
\forall w \in \tilde{V}_h^k \; \displaystyle \sup_{v \in V_h^h \setminus \{0\}} \frac{a_h(w,v)}{\| v \|_{1,h}} \geq \alpha \| w \|_{\tilde{W},2,h}
\end{array}
\right.
\end{equation}
\end{theorem}
\prov \eqref{infsup} is a mere consequence of \eqref{coercivbarah}, \eqref{coercivah} and \eqref{equiv}, with $\bar{\alpha} = C^{-} \bar{\alpha}^{'}$ and $\alpha = C^{-} \alpha^{'}$. \QED

\subsection{Uniform boundedness}

We pursue our reliability analysis by addressing the uniform continuity of the bilinear form $\bar{a}_h$. We have
\begin{theorem} 
\label{boundbarah}
There exists a constant $\bar{A}$ independent of $h$ such that 
\begin{equation}
\label{continuitybarah}
\bar{a}_h(w,v) \leq \bar{A} \| w \|_{\tilde{W},2,h} \| v \|_{1,h} \; \forall w \in \tilde{W}^2_h \mbox{ and } \forall v \in V^k_h.
\end{equation}
\end{theorem}
\prov
Let us first consider the case of $\bar{b}_h$. We have 
\begin{equation}
\label{boundbarbh1}
\bar{b}_h(w,v) \leq \displaystyle \max[\| \tilde{d} \|_{0,\infty,\Omega_h},[\| \tilde{r} \|_{0,\infty,\Omega_h}] \| w \|_{1,h} \| v \|_{1,h} 
+ \| \eta \|_{0,\infty,\Gamma} \displaystyle \sum_{T \in {\mathcal S}_h} \| v \|_{0,e_T} \| w \|_{0,e_T}.
\end{equation}
Thus, it turns out that
\begin{equation}
\label{boundbarbh2}
\bar{b}_h(w,v) \leq \displaystyle \max[\| \tilde{d} \|_{0,\infty,\tilde{\Omega}},\| \tilde{r} \|_{0,\infty,\tilde{\Omega}},\| \eta \|_{0,\infty,\Gamma} ] [\| w \|_{1,h} \| v \|_{1,h} +  \| v \|_{0,\Gamma_h} \| w \|_{0,\Gamma_h}].
\end{equation}
Taking into account the Trace Theorem for $H^1(\Omega_h)$, the following bound derives from \eqref{boundbarbh2}, where the constant 
$\bar{C}_b$ depends on $\tilde{d}$, $\tilde{r}$ and $\eta$ but not on $h$.
\begin{equation}
\label{boundbarbh}
\bar{b}_h(w,v) \leq \bar{C}_b \| w \|_{1,h}  \| v \|_{1,h} \; \forall w \in  \tilde{W}^2_h \mbox{ and } \forall v \in V^k_h. 
\end{equation}
Next we endeavor to find an upper bound for $\bar{c}_h(w,v)$ other than \eqref{estimbarch} in order to accommodate functions $w$ in $H^2(\tilde{\Omega}_h)$.
Recalling the expression \eqref{barch} of $\bar{c}_h$ we first note that it extends to $w \in H^2(\tilde{\Omega}_h)$, in which case we have
\begin{equation}
\label{barch1}
\left\{
\begin{array}{l}
\displaystyle \int_{e_T} {\mathcal F}_T(w,v) \leq  \displaystyle \left\{ \| \eta \|_{0,\infty,\Gamma} \displaystyle \sqrt{\int_{e_T} \left( \int_M^P   \left|\frac{\partial w}{\partial n_h} \right| \right)^2}  +  \| \tilde{d} \|_{0,\infty\tilde{\Omega}} \displaystyle  \left[
\sqrt{\int_{e_T}   \left( \int_M^P   \left|\frac{\partial^2 w}{\partial n_P \partial n_h}\right|\right)^2} \right. \right. \\
\left. \left. + \displaystyle s_m h_T \left\| \frac{\partial w}{\partial \tau_P} \right\|_{0,e_T} \right] + 
  | \tilde{d} |_{1,\infty\tilde{\Omega}} \displaystyle 
\sqrt{\int_{e_T}   \left( \int_M^P   \left|\frac{\partial w}{\partial n_P}\right| \right)^2}
\right\}\| v \|_{0,e_T}.  
\end{array}
\right.
\end{equation}
Recalling \eqref{partialtauP} and \eqref{partialnP}, 
after straightforward calculations it follows from \eqref{barch1} that
\begin{equation}
\label{barch2}
\left\{
\begin{array}{l}
\int_{e_T} {\mathcal F}_T(w,v) \leq \left\{\max[\|\eta \|_{0,\infty,\Gamma},\|\tilde{d} \|_{1,\infty,\tilde{\Omega}}] \displaystyle \max[  C_{\Gamma}^{1/2}, s_m] h_T \right.\\
\left. \left( \| \nabla w \|_{1,\Delta_T} + s_m h_T \displaystyle \left\| \frac{\partial w}{\partial n_h} \right\|_{0,e_T} + 
 \left\| \frac{\partial w}{\partial \tau_h} \right\|_{0,e_T} \right)  \right\}  \| v \|_{0,e_T}. 
\end{array}
\right.
\end{equation}
Sweeping ${\mathcal S}_h$ and denoting by $h_m$ the maximum value of $h$ in ${\mathcal P}$, from \eqref{barch2} and \eqref{barch} we easily obtain 
\begin{equation}
\label{barch3}
\bar{c}_h(w,v) \leq \max[\| \eta \|_{0,\infty,\Gamma}, \| \tilde{d} \|_{1,\infty,\tilde{\Omega}}] \max[C^{1/2}_{\Gamma}\!,s_m]   
\displaystyle \left(\! \| w \|_{\tilde{W},2,h} \!+\! h \!\sqrt{1+s_m^2h_m^2} \| \nabla w \|_{0,\Gamma_h}\! \right) \| v \|_{0,\Gamma_h}.
\end{equation}
Resorting to the Trace Theorem, we can assert that there is a constant $C_{tr}^{'}$ supposedly (cf. Remark 1) independent of $h$  such that
$$\| \nabla w \|_{0,\Gamma_h} \leq C_{tr}^{'} \| w \|_{2,\Omega_h}.$$ 
This readily yields for a mesh-independent constant $\bar{C}(d,\eta)$
\begin{equation}
\label{barch4}
\bar{c}_h(w,v) \leq \bar{C}(\eta,\tilde{d}) \| w \|_{\tilde{W},2,h} \| v \|_{1,h} \; \forall w \in H^2(\Omega_h) \mbox{ and } \forall v \in V^k_h.
\end{equation}
Finally, recalling \eqref{estimbarch}, we obtain
\begin{equation}
\label{estimabarch}
\bar{c}_h(w,v) \leq \tilde{C}_c \| w \|_{\tilde{W},2,h} \| v \|_{1,h} \; \forall (w;v) \in  \tilde{W}^2_h \times V^k_h,  
\end{equation}
where $\tilde{C}_c = \max[\bar{C}_c h_m^{1/2},\bar{C}(\eta,\tilde{d})]$.
Combining \eqref{boundbarbh} and \eqref{estimabarch}, \eqref{continuitybarah} is seen to hold $\forall w \in \tilde{W}^2_h \times V_h^k$ with $\bar{A} = \bar{C}_b + \tilde{C}_c$. 
\QED

\subsection{Convergence results}

On the basis of the properties of bilinear forms $\bar{a}_h$ and $a_h$ established in the previous subsections, it is possible to demonstrate the convergence of the approximation method \eqref{vph} of the model problem \eqref{eq}, as long as the solution of the latter is sufficiently smooth.\\
With this aim we first note that the fictive solution method \eqref{barvph} was designed to be conforming. More precisely, we have
\begin{e-proposition}
\label{residual} 
The variational residual $R(w,v,\tilde{f},g):=\bar{a}_h(w,v) - \bar{L}_h(v)$ vanishes identically for $w = \tilde{u}$.\\ 
\end{e-proposition}

\prov Integrating by parts in $\Omega_h$, since by construction $-\nabla \cdot \tilde{d} \nabla \tilde{u} + \tilde{r} \tilde{u} - \tilde{f} \equiv 0$ in $\Omega_h$ we obtain 
\begin{equation}
\label{residual0}
\left\{
\begin{array}{l}    
R(\tilde{u},v,\tilde{f},g) = \displaystyle \sum_{T \in {\mathcal S}_h} \int_{e_T} \left\{ \tilde{d} \frac{\partial \tilde{u}}{\partial n_h}  \right.  \\
\left.  + \mbox{ cos } \theta_T \displaystyle \left[\bar{\eta} \overline{\tilde{u}} + \; \displaystyle \int_M^P \frac{\partial}{\partial n_h}\left(\tilde{d} \frac{\partial \tilde{u}}{\partial n_P}\right) - \bar{g} \right] + \displaystyle \tilde{d} \mbox{ sin } \theta_T \frac{\partial \tilde{u}}{\partial \tau_P} \right\} v.  
\end{array}
\right.
\end{equation}
Now, after straightforward calculations, we verify that 
\begin{equation}
\label{residual1} 
\tilde{d} \frac{\partial \tilde{u}}{\partial n_h} + \mbox{ cos } \theta_T \int_M^P \frac{\partial}{\partial n_h}\left(\tilde{d} \frac{\partial \tilde{u}}{\partial n_P}\right) + \tilde{d} \displaystyle \mbox{ sin } \theta_T \frac{\partial \tilde{u}}{\partial \tau_P} = \displaystyle  \bar{d} \mbox{ cos } \theta_T \overline{\frac{\partial \tilde{u}}{\partial n_P}} \mbox{ on } e_T. 
\end{equation}
Plugging \eqref{residual1} into \eqref{residual0} we come up with  
\begin{equation}
\label{residual2}
R(\tilde{u},v,\tilde{f},g) = \displaystyle \sum_{T \in {\mathcal S}_h} \int_{e_T} \mbox{ cos } \theta_T \left(\bar{\eta} \overline{\tilde{u}} + \displaystyle \bar{d} \overline{\frac{\partial \tilde{u}}{\partial n}} - \bar{g} \right)  v.  
\end{equation}     
On the other hand, owing to the boundary conditions in \eqref{eq}, we know that
\begin{equation}
\label{residual3}
\eta \tilde{u} + \displaystyle d \frac{\partial \tilde{u}}{\partial n} - g = \eta u + \displaystyle d \frac{\partial u}{\partial n} - g = 0 \mbox{ on } \Gamma_T \; \forall T \in {\mathcal S}_h.
\end{equation} 
Since the transformation of the left hand side of \eqref{residual3} onto $e_T$ necessarily vanishes identically as well, $R(\tilde{u},v,\tilde{f},g) = 0 \; \forall v \in V^k_h$. 
\QED \\    

Using Proposition \ref{residual}, Theorem 5.1 in \cite{COAM} (cf. equation (20)) allows us to write 
\begin{equation}
\label{estimbarvph0}
\| \bar{u}^k_h - \tilde{u} \|_{\tilde{W},2,h} \leq \displaystyle \frac{\bar{A}}{\bar{\alpha}} \| \tilde{u} - \tilde{\pi}_h^k(\tilde{u}) \|_{\tilde{W},2,h},
\end{equation}
where $\tilde{\pi}_h^k(w)$ is the standard $\tilde{V}_h^k$-interpolate of a function $w \in H^2(\tilde{\Omega})$.\\
On the basis of the classical interpolation theory applying also to triangles with a curved edge (cf. \cite{BrennerScott}) we easily infer that, as long as $\tilde{u} \in H^{k+1}(\tilde{\Omega})$, for a constant $\bar{C}_k$ independent of $h$ we have  
\begin{equation}
\| \tilde{u} - \tilde{\pi}_h^k(\tilde{u}) \|_{\tilde{W},2,h} \leq \bar{C}_k h^k \| \tilde{u} \|_{k+1,\tilde{\Omega}}.
\end{equation}
In short, we have proved 
\begin{theorem}
\label{barerrestim}
Il holds 
\begin{equation}
\label{estimbarvph}
\| \bar{u}^k_h - \tilde{u} \|_{\tilde{W},2,h} \leq \displaystyle \frac{\bar{A}\bar{C}_k}{\bar{\alpha}} h^k \| \tilde{u} \|_{k+1,\tilde{\Omega}}.
\end{equation}
\QED
\end{theorem}
Further, exploiting the relation $\| \tilde{u} - \bar{u}_h^k \|_{m,\tilde{T}} \leq \| \tilde{u} - \tilde{\pi}_h^k(\tilde{u}) \|_{m,\tilde{T}} + \| \tilde{\pi}_h^k(\tilde{u}) - \bar{u}_h^k \|_{m,\tilde{T}}$ for $T \in {\mathcal T}_h$ with $1 \leq m \leq k$, together with the inverse inequalities $\| w \|_{m,\tilde{T}} \leq \tilde{C}_{I,m} h_T^{1-m} \| w \|_{1,\tilde{T}}$ with a constant $\tilde{C}_{I,m}$ independent of $T$ for any $w \in H^m(\tilde{T})$, without any difficulty we come up with the following estimates for constants $\tilde{C}_{k,m}$ and $\tilde{C}_{m}$ independent of $h$ and $\tilde{u}$.
\begin{equation}
\label{estimatemk}
\left\{
\begin{array}{l}
\| \tilde{u} - \bar{u}_h^k \|_{m,\tilde{\Omega}_h} \leq \tilde{C}_{k,m} h^{k-m+1} \| \tilde{u} \|_{k+1,\tilde{\Omega}} \\
\\
\| \bar{u}_h^k \|_{m,\tilde{\Omega}_h} \leq \tilde{C}_{m} \| \tilde{u} \|_{m,\tilde{\Omega}_h} \mbox{ for } m=1,\ldots,k.
\end{array}
\right.
\end{equation} 
Now, since the existence and uniqueness of $u_h^k$ is ensured by \eqref{infsup} we are ready to establish an error estimate for it. With this aim it is appropriate to apply the error bound (30) given in \cite{COAM} taking into account \eqref{estimbarvph}. In accordancence with both inequalities we can write,
\begin{equation}
\label{errorbound0}
\| \tilde{u} - u^k_h \|_{\tilde{W},2,h} \leq \displaystyle \frac{\bar{A}\bar{C}_k}{\bar{\alpha}} h^k \| \tilde{u} \|_{k+1,\tilde{\Omega}} + \displaystyle \frac{1}{\alpha} \sup_{v \in V_h^k \setminus \{0\}} \frac{|a_h(\bar{u}^k_h,v)-\bar{a}_h(\bar{u}^k_h,v)|+|L_h(v)-\bar{L}_h(v)|}{\| v \|_{1,h}}.
\end{equation}
\eqref{errorbound0} tells us that the estimation of the difference between the bilinear forms $\bar{a}_h$ and $a_h$ applied to $\bar{u}_h^k$ and $v$ and the linear forms $\bar{L}_h$ and $L_h$ applied to $v \in V_h^k$ is all that remains to be done.\\
First of all we have
\begin{e-proposition}
\label{rhs}
Provided $f \in H^k(\Omega)$ $g \in H^k(\Gamma)$ and $\Gamma$ is of the $C^{k+1}$-class, there exists a constant $C_L^{k}(\tilde{f},g)$ independent of $h$ such that
\begin{equation}
\label{barLh-Lh}
| \bar{L}_h(v) - L_h(v) | \leq C_L^{k}(\tilde{f},g) h^k \| v \|_{1,h} \; \forall v \in V_h^k.
\end{equation} 
\end{e-proposition}
\prov
By the standard interpolation theory it is clear that 
\begin{equation}
\label{errorf}
\int_T [\tilde{f}-\rho_{k-1}(\tilde{f})]v \leq C_{\rho}h_T^k | \tilde{f} |_{k,T} \| v \|_{0,T} \; \forall T \in {\mathcal T}_h, 
\end{equation}
where $C_{\rho}$ depends neither on $T$ nor on $f$ and $v$.\\
On the other hand we also have for $\phi= \mbox{ cos } \theta_T \bar{g}$
\begin{equation}
\label{errorg0}
\int_{e_T} [\phi-\sigma_{k-1}(\phi)]v \leq C_{\sigma} h_T^k | \phi |_{k,e_T} \| v \|_{0,e_T} \; \forall T \in {\mathcal S}_h, 
\end{equation}
Now we observe that, on the basis of Proposition 2.2 of \cite{ZAMM} and the chain rule, we can assert that for $l=0,1,\ldots,k$ there exist two constants $C_l(\Gamma)$ and $\bar{C}_l(\Gamma)$ such that $| \mbox{cos }\theta_T|_{l,\infty,e_T} \leq C_l(\Gamma)$ and $| \bar{g} |_{l,e_T} \leq \bar{C}_l(\Gamma) \| g \|_{l,\Gamma_T}$. It easily follows that $| \phi |_{k,e_T} \leq C_k(\Gamma) \| g \|_{k,\Gamma_T}$ for a suitable constant $C_k(\Gamma)$ independent of $T$. In short, it holds  
 \begin{equation}
\label{errorg}
\int_{e_T} [\phi-\sigma_{k-1}(\phi)]v \leq C_{\sigma} {\mathcal C}_k h_T^k \| g \|_{k,\Gamma_T} \| v \|_{0,e_T} \; \forall T \in {\mathcal S}_h \mbox{ with } \phi = \mbox{ cos } \theta_T \bar{g}.
\end{equation} 
By summation over ${\mathcal S}_h$, the Cauchy-Schwarz inequality and the Trace Theorem, \eqref{barLh-Lh} easily derives from \eqref{errorf} and \eqref{errorg}. \QED \\  

\begin{e-proposition}
\label{lhs}
For $k \geq 2$, assume that $\tilde{d} \in C^{2k-1}(\bar{\Omega})$, $\tilde{r} \in C^{2k-1}(\bar{\Omega})$, $\eta \in C^{2k-1}(\Gamma)$ and $\Gamma$ is of the $C^{2k}$-class. If the numerical quadrature formulae ${\mathcal J}^k_T$ and ${\mathcal I}^k_T$ integrate exactly polynomials of degree respectively less than or equal to $2k-2$ in every $T \in {\mathcal T}_h$ and of degree less than or equal to $2k-1$ in $e_T$ for every $T \in {\mathcal S}_h$, then there exists a constant ${\mathcal C}_A$ depending on $\tilde{d}$, $\tilde{r}$, $\eta$ and $\Gamma$ but not on $h$, such that
\begin{equation}
\label{ahbarah}
| \bar{a}_h(w,v) - a_h(w,v) | \leq  {\mathcal C}_A h^k \| w \|_{k,h} \| v \|_{1,h} \; \forall w \in \tilde{V}_h^k \mbox{ and }\forall v \in V_h^k.\\
\end{equation} 
\end{e-proposition}

\prov
First we consider the case of $b_h$ and $\bar{b}_h$. Akin to previous results given in this work using the master element $\hat{T}$, we can 
assert that there exists a constant $\bar{C}_{B}$ independent of $h$ such that 
\begin{equation} 
\label{bhbarbh0}
\left\{
\begin{array}{l}
|\bar{b}_h(w,v) - b_h(w,v)|  \leq C_{B} \displaystyle \left[ \sum_{T \in {\mathcal T}_h} h_T^{2k-1} \left(| \tilde{d} \nabla w \cdot \nabla v |_{2k-1,1,T} + | \tilde{r}  w  v |_{2k-1,1,T} \right) \right. \\
\left. +  \displaystyle \sum_{T \in {\mathcal S}_h} |{\mathcal H}_T^{2k-1}(\mbox{cos }\theta_T \bar{\eta} w v)| \right].
\end{array}
\right.
\end{equation}
It is clear that there are two constants $\tilde{C}^{'}_{D}(\tilde{d})$ and $\tilde{C}^{'}_{R}(\tilde{r})$ not depending on $T$, such that
\begin{equation}
\label{bhbarbh1}
\left\{
\begin{array}{l}
 | \tilde{d} \nabla w \cdot \nabla v |_{2k-1,1,T} \leq \tilde{C}^{'}_{D}(\tilde{d}) \displaystyle \sum_{l=0}^{k-1} \sum_{j=0}^{l}  
\| w \|_{j+1,T} \| v \|_{l-j+1,T}  \\
 \mbox{and} \\
 | \tilde{r}  w v |_{2k-1,1,T} \leq \tilde{C}^{'}_{R}(\tilde{r}) \displaystyle \left( \| w \|_{k,T} \| v \|_{0,T} +
 \sum_{l=1}^{k-1} \sum_{j=0}^{l} \| w \|_{j,T} \| v \|_{l-j+1,T} \right).   
 \end{array}
\right.
\end{equation}
Using classical inverse inequalities we further come up with two other constants $C^{''}_{D}(\tilde{d})$ and $C^{''}_{R}(\tilde{r})$ such that  
\begin{equation}
\label{bhbarbh1bis}
\left\{
\begin{array}{l}
| \tilde{d} \nabla w \cdot \nabla v |_{2k-1,1,T} \leq C^{''}_{D}(\tilde{d}) \displaystyle \sum_{l=0}^{k-1} \sum_{j=0}^{l}
\| w \|_{j+1,T} h_T^{-l+j} \| v \|_{1,T}  \\
 \mbox{and} \\
| \tilde{r}  w v |_{2k-1,1,T} \leq C^{''}_{R}(\tilde{r}) \displaystyle \left( \| w \|_{k,T} \| v \|_{0,T} +
 \displaystyle \sum_{l=1}^{k-1} \sum_{j=0}^{l} \| w \|_{j,T} h_T^{-l+j} \| v \|_{1,T} \right).   
\end{array}
\right.
\end{equation}
or yet for constants $C_{D}(\tilde{d})$ and $C_{R}(\tilde{r})$ of the same nature
\begin{equation}
\label{bhbarbh1ter}
\left\{
\begin{array}{l}
h_T^{2k-1} | \tilde{d} \nabla w \cdot \nabla v |_{2k-1,1,T} \leq C_{D}(\tilde{d}) h_T^{k} \| w \|_{k,T} \| v \|_{1,T}  \\
 \mbox{and} \\
h_T^{2k-1}| \tilde{r}  w v |_{2k-1,1,T} \leq C_{R}(\tilde{r}) h_T^{k} \| w \|_{k,T} \| v \|_{1,T}.   
 \end{array}
\right.
\end{equation}
Next, we resort to Lemma \ref{auxiliary4} in order to estimate the term ${\mathcal H}^{2k-1}_T(\phi w v)$ with $\phi = \mbox{cos } \theta_T \bar{\eta}$ and $z= w$. Using repeatedly elementary Calculus, we establish right away the uniform bound \eqref{d2kphidtau2k} given below.\\
Indeed, recalling the arguments in Proposition 2.2 of \cite{ZAMM}, first we observe that $| \mbox{cos } \theta_T |_{l,\infty,e_T}$ is uniformly bounded above by a constant $C_{2k}^{'}(\Gamma)$, from $l=0$ through $l=2k$. On the other hand, we remind that $\mbox{d}\bar{\eta}/\mbox{d}\tau_h=\mbox{ cos } \theta_T \overline{\mbox{d}\eta/\mbox{d}\tau}$ along $e_T$. Hence, making repeated use of the above defined constant $C_{2k}^{'}(\Gamma)$, after straightforward calculations, we come up with another constant $\bar{C}_{2k}(\Gamma)$ depending only on $\Gamma$ such that $\| \bar{\eta} \|_{2k,\infty,e_T} \leq \bar{C}_{2k}(\Gamma) \| \eta \|_{2k,\infty,\Gamma}$.\\
Summarizing, setting $C_{2k}(\Gamma)= {\bf C}_{2k}^k C^{'}_{2k}(\Gamma) \bar{C}_{2k}(\Gamma)$, where ${\bf C}_n^l$ represents the integer $n \; choose \; l$ for $l \leq n$, we can assert that
\begin{equation}
\label{d2kphidtau2k}
\| \phi \|_{2k,\infty,e_T} = {\bf C}_{2k}^k \displaystyle \sum_{l=0}^{2k} | \mbox{cos } \theta_T |_{l,\infty,e_T} | \eta |_{2k-l,\infty,e_T} 
\leq C_{2k}(\Gamma) \| \eta \|_{2k,\infty,e_T}. 
\end{equation}
It follows that there exists a constant $C^{'}_H(\eta)$ expressed only in terms of $\| \eta \|_{2k,\infty,\Gamma}$ and $\Gamma$ itself, such that
\begin{equation}
\label{bhbarbh2}
{\mathcal H}^{2k-1}_T(\mbox{cos } \theta_T \bar{\eta} w v) \leq C^{'}_H(\eta) h^{2k} \displaystyle \left( \| w \|_{k,e_T} \| v \|_{0,e_T} + \sum_{l=1}^{k-1} \sum_{j=0}^l \| w \|_{j,e_T} \| v \|_{l-j+1,e_T} \right).
\end{equation}
Now, using inverse inequalities for fractional Sobolev spaces given in \cite{Georgoulis}, we obtain for a constant $C^{''}_H(\eta)$  
\begin{equation}
\label{bhbarbh2bis}
{\mathcal H}^{2k-1}_T(\mbox{cos } \theta_T \bar{\eta} w v) \leq C^{''}_H(\eta) h_T^{2k} \displaystyle \left( \| w \|_{k,e_T} \| v \|_{0,e_T} + \sum_{l=1}^{k-1} \sum_{j=0}^l \| w \|_{j,e_T} h_T^{j-l-1/2} \| v \|_{1/2,e_T} \right),
\end{equation}
or yet, for another similar constant $C_H(\eta)$
\begin{equation}
\label{bhbarbh2ter}
{\mathcal H}^{2k-1}_T(\mbox{cos } \theta_T \bar{\eta} w v) \leq C_H(\eta) h_T^{k+1/2} \| w \|_{k,e_T} \| v \|_{1/2,e_T}.
\end{equation}
Putting together \eqref{bhbarbh0}, \eqref{bhbarbh1ter} and \eqref{bhbarbh2ter} we promptly obtain for a constant $C_B$ depending on $\tilde{d}$, $\tilde{r}$ and $\eta$, but not on $h$  
\begin{equation} 
\label{bhbarbh3}
|\bar{b}_h(w,v)\! -\! b_h(w,v)| \leq C_B \displaystyle \left( \sum_{T \in {\mathcal T}_h} h_T^{k} \| w \|_{k,T} \| v \|_{1,T} +  \displaystyle \sum_{T \in {\mathcal S}_h} h_T^{k+1/2} \| w \|_{k,e_T} \| v \|_{1/2,e_T} \right).
\end{equation}
Plugging into \eqref{bhbarbh3} the inequalities 
\begin{equation}
\label{inverse} 
\| p \|_{k,e_T} \leq h_T^{1/2} \| p \|_{k,\infty,T} \leq C_I^k h_T^{-1/2} \| p \|_{k,T} \; \forall p \in P_k(T),
\end{equation}
for a constant $C_I^k$ not depending on $T$, owing to the Trace Theorem for the space $H^1(\Omega_h)$, we readily come up with a mesh-independent constant ${\mathcal C}_B$ such that
\begin{equation} 
\label{bhbarbh4}
|\bar{b}_h(w,v) - b_h(w,v)| \leq {\mathcal C}_B h^{k} \| w \|_{k,h} \| v \|_{1,h}. 
\end{equation}
\indent Next we address the case of $\bar{c}_h$ and $c_h$. We have
\begin{equation}
\label{chbarch0}
\left\{
\begin{array}{l}
|\bar{c}_h(w,v) - c_h(w,v) | \leq \displaystyle \sum_{T \in {\mathcal S}_h} \left[ \bar{c}_T(w,v) - c_T(w,v) \right].\\
\mbox{where } \bar{c}_T \mbox{ is defined in } \eqref{barchwv1} \mbox{ and}\\
c_T(w,v)= \displaystyle {\mathcal I}^k_T \left\{ \left[ \mbox{cos } \theta_T \left( \bar{\eta}\left( \bar{w}- w \right) \; + \; \displaystyle \int_M^P \frac{\partial}{\partial n_h}\left(\tilde{d} \frac{\partial w}{\partial n_P}\right) \right) \; 
+ \; \tilde{d} \mbox{ sin } \theta_T \displaystyle \frac{\partial w}{\partial \tau_P} \right] v \right\}.
\end{array}
\right.
\end{equation}
Here it is convenient to rewrite the integral along the segment $\overline{MP}$ orthogonal to $e_T$ for $M \in e_T$ and $P \in \Gamma_T$ in the expressions of both $c_T$ and $\bar{c}_T$ as 
$$\displaystyle \left[\bar{d} \overline{\frac{\partial w}{\partial n_P}}\right](M)-\left[\tilde{d} \frac{\partial w}{\partial n_P}\right](M).$$ 
In doing so we can handle all the resulting terms in the very same manner as $\mbox{cos } \theta_T \bar{\eta}w$ on $e_T$ in equation 
\eqref{bhbarbh2}. As a matter of fact, in the light of \eqref{dudnh1} we can simplify a little the expressions of $\bar{c}_T$ and $c_T$, thereby obtaining
\begin{equation}
\label{chbarch1}
\left\{
\begin{array}{l}
\bar{c}_T(w,v) -  c_T(w,v) = {\mathcal H}_T^{2k-1}( \bar{\mathcal G}_T(w) v)+{\mathcal H}_T^{2k-1}({\mathcal G}_T(w) v ) \\
\mbox{with} \\
\bar{\mathcal G}_T(w):= \mbox{cos } \theta_T \displaystyle \left(\bar{\eta} \bar{w} \; + \; \bar{d} \displaystyle \overline{\frac{\partial w}{\partial n}} \right) \\
\mbox{and} \\ 
{\mathcal G}_T(w):= - \mbox{cos } \theta_T \bar{\eta} w \; - \; \tilde{d} \displaystyle \frac{\partial w}{\partial n_h}.
\end{array}
\right.
\end{equation}
Mimicking the estimate \eqref{bhbarbh2ter} we infer the existence of a constant $C_{G}(\tilde{d})$ independent of $T$ such that 
\begin{equation}
\label{chbarch2}
{\mathcal H}_T^{2k-1}({\mathcal G}_T(w) v ) \leq  h_T^{k+1/2} \| v \|_{1/2,e_T}\displaystyle \left( C_H(\eta)  \| w \|_{k,e_T}  + 
C_D(\tilde{d}) \displaystyle \left\| \frac{\partial w}{\partial n_h}\right\|_{k-1,e_T}  \right).
\end{equation}
Using the inequalities \eqref{inverse}, the equation given next trivially derives from \eqref{chbarch2}. 
\begin{equation}
\label{chbarch3}
{\mathcal H}_T^{2k-1}({\mathcal G}_T(w) v ) \leq  C_I^k h_T^{k} \| v \|_{1/2,e_T}\displaystyle \left( C_H(\eta)  \| w \|_{k,T}  + 
C_D(\tilde{d}) \| \nabla w \|_{k-1,T}  \right).
\end{equation}
The case of $ {\mathcal H}_T^{2k-1}( \bar{\mathcal G}_T(w) v)$ is more complex, since neither $\bar{w}$ nor $\overline{\partial w/\partial n_P}$ are polynomials. Nevertheless, here again we can resort to Proposition 2.2 of \cite{ZAMM} and for better insight we split $\bar{\mathcal G}_T(w)$ into the sum of $\bar{\mathcal G}^{a}_T(w)$ and $\bar{\mathcal G}^{b}_T(w)$ given by.
\begin{equation}
\label{splitbarG}
\left\{
\begin{array}{l}
\bar{\mathcal G}^{a}_T(w):= \mbox{cos } \theta_T \bar{\eta} \bar{w} \\
\mbox{and} \\ 
\bar{\mathcal G}^{b}_T(w):= \mbox{cos } \theta_T \; \bar{d} \displaystyle \overline{\frac{\partial w}{\partial n}}.
\end{array}
\right.
\end{equation}
In order to apply Lemma \ref{auxiliary4} to the term ${\mathcal H}_T^{2k-1}(\bar{\mathcal G}^{a}_T(w)v)$, like in the case of 
the first term of ${\mathcal G}_T(w)$ in \eqref{chbarch2}, it is also convenient to take here $\phi = \bar{\eta} \mbox{ cos } \theta_T$. In doing so, now we have $z=\bar{w}$ and thus, analogously to \eqref{bhbarbh2ter}, it holds 
\begin{equation}
\label{chbarch4}
{\mathcal H}_T^{2k-1}(\bar{\mathcal G}^{a}_T(w) v) \leq C_H(\eta) h_T^{k+1/2} \| z \|_{2k,e_T}  \| v \|_{1/2,e_T} .
\end{equation} 
Akin to \eqref{dtauh}, in order to estimate the norms $\|z\|_{l,e_T}$ for $l=1,\ldots,2k$, we must handle a compound function of the form $\bar{\varphi}(x)=\varphi(x,y(x))$, for suitable $2k$-times differentiable functions $\varphi$ and $y$, $x=\tau_h$, $(x,0)$ and $(x,y(x))$ being the local coordinates in a direct cartesian frame of $M$ and $P$, respectively. \\
In the estimation of the derivatives $\mbox{d}^lz/\mbox{d}\tau_h^l$ for $l=1,\ldots,2k$ that follows, with $z = \bar{w}$, 
we use the notation $\varphi^{'}$ and $\varphi^{''}$ for the first and second derivatives of a function $\varphi$ of the variable $x$. In doing so, the total first order derivative of $z$ along $e_T$ is expressed by 
\begin{equation}
\label{derivatotal10}
[\mbox{d}z/\mbox{d}\tau_h](x) = [\partial w/\partial \tau_h](x,y(x)) + [\partial w/\partial n_h](x,y(x))y^{'}(x).
\end{equation} 
As already pointed out in \eqref{dtauh}, $\theta_T = -\mbox{atan }y^{'}$ and furthermore $\mbox{cos }\theta_T = [1+(y^{'})^2]^{-1/2}$. It easily follows that 
\begin{equation}
\label{derivatotal1}
\left\{
\begin{array}{l}
\displaystyle \left[\frac{\mbox{d}z}{\mbox{d}\tau_h}\right](M)= a_{11}(M) \displaystyle \left[\frac{\partial w}{\partial \tau}\right](P), \\
\mbox{where}\\
a_{11}(M) = \displaystyle \left[\sqrt{1+ \left(y^{'}\right)^2}\right](M).
\end{array}
\right.
\end{equation}
From the regularity assumption on $\Gamma$ by a straightforward calculation, it is easily seen that $a_{11} \in W^{2k-1,\infty}(e_T)$ and moreover $\| a_{11} \|_{2k-1,\infty,e_T}$ is bounded above by a constant $C_{11}(\Gamma)$ depending only on $\Gamma$.\\ 
Now, switching to the second order derivative of $z$ along $e_T$, without any difficulty we obtain 
\begin{equation}
\label{derivatotal20}
\left\{
\begin{array}{l}
\displaystyle \left[\frac{\mbox{d}^2z}{\mbox{d}\tau_h^2}\right](M) = \displaystyle \left\{[\mbox{tan }\theta_T](M) \displaystyle 
\left[\frac{y^{''}}{1+\left(y^{'}\right)^2} \right](M) \displaystyle \left[\frac{\partial w}{\partial \tau}\right](P) \right. \\
\left. + \displaystyle \left[\frac{\partial^2 w}{\partial \tau \partial \tau_h}\right](P) + \left[\frac{\partial^2 w}{\partial \tau \partial n_h} \right](P) \displaystyle \left[ y^{'} \right](M) \right\} \left[\mbox{cos}^{-1}\theta_T \right](M).
\end{array}
\right. 
\end{equation}
After straightforward calculations we come up with
\begin{equation}
\label{derivatotal2}
\left\{
\begin{array}{l}
\displaystyle \left[\frac{\mbox{d}^2z}{\mbox{d}\tau_h^2}\right](M) = a_{21}(M) \displaystyle \left[\frac{\partial w}{\partial \tau}\right](P) + a_{22}(M) \displaystyle \left[\frac{\partial^2 w}{\partial \tau^2}\right](P), \\
\mbox{where} \\
a_{21}(M) = \displaystyle \left[y^{'} \kappa\right](M) \mbox{ and } a_{22}(M) = \displaystyle \left[1+ \left(y^{'}\right)^2 \right](M), \\
\end{array}
\right. 
\end{equation}
$\kappa(M) = \left\{\displaystyle\frac{y^{''}}{\left[1+\left(y^{'}\right)^2\right]^{3/2}} \right\}(M)$ being the curvature of $\Gamma$ at the point $P$ with local coordinates $(x,y(x))$.\\
Akin to $a_{11}$, by the assumed regularity of $\Gamma$ it is not difficult to establish that both $a_{21}$ and $a_{22}$ belong to $W^{2k-2,\infty}(e_T)$ and moreover there are two constants $C_{21}(\Gamma)$ and $C_{22}(\Gamma)$ depending only on $\Gamma$ such that 
$\| a_{2i} \|_{2k-2,\infty,e_T} \leq C_{2i}(\Gamma)$ for $i=1,2$.\\ 
Next we proceed by mathematical induction: For a given integer $l$ comprised between $3$ and $2k-1$, assume that there exist functions $a_{jl} \in W^{2k-l,\infty}(e_T)$ and constants $C_{jl}(\Gamma)$ depending only on $\Gamma$, for $j=1,2,\ldots,l$, such that 
\begin{equation}
\label{dlzdtauhl}
\left\{
\begin{array}{l}
\displaystyle \frac{\mbox{d}^lz}{\mbox{d}\tau_h^l} = \displaystyle \sum_{j=1}^l a_{jl} \frac{\partial^j w}{\partial \tau^j} \\
\mbox{with} \\
\| a_{jl} \|_{2k-l,\infty,e_T} \leq C_{jl}(\Gamma) \mbox{ for } j=1,2,\ldots,l.
\end{array}
\right.
\end{equation}
We have to prove that two relations of the same type as \eqref{dlzdtauhl} hold true for the integer $l+1 \leq 2k$. With this aim we combine  
\eqref{derivatotal1} and \eqref{dlzdtauhl} 
\begin{equation}
\label{dl1zdtauhl1}
\left\{
\begin{array}{l}
\displaystyle \left[\frac{\mbox{d}^{l+1}z}{\mbox{d}\tau_h^{l+1}}\right](M)= 
\displaystyle \left\{ \frac{\mbox{d} \displaystyle \left[\frac{\mbox{d}^lz}{\mbox{d}\tau_h^l}\right]}{\mbox{d}\tau_h}\right\}(M)\\
=  \displaystyle \sum_{j=1}^l \left[ a^{'}_{jl}(M) \left(\frac{\partial^j w}{\partial \tau^j}\right)(P) + a_{11}(M)a_{jl}(M)
\left(\frac{\partial^{j+1} w}{\partial \tau^{j+1}}\right)(P) \right].
\end{array}
\right.
\end{equation} 
This leads to \eqref{derivatotal1} with $a_{1,l+1} = a^{'}_{1l}$, $a_{j,l+1} = a^{'}_{jl}+a_{11}a_{j-1,l}$ for $j=2,\ldots,l$ and 
$a_{l+1,l+1}=a_{11}a_{ll}$. It is easy to figure out that, for $j=1,2,\ldots,l+1$, there exist constants $C_{j,l+1}(\Gamma)$ depending only on $\Gamma$ such that $\| a_{j,l+1} \|_{2k-l-1,\infty,e_T} \leq C_{j,l+1}(\Gamma)$. It follows that \eqref{dlzdtauhl} holds for any $l$ between $2$ and $2k$.\\
Now, for a given $l \geq 0$, we denote by $D^lz$ the $l$-th order tensor, whose components are the $2^l$ $l$-th order partial derivatives of a  function $z \in C^{2k}(\overline{\tilde{\Omega}})$. Noticing that for $j \geq 1$, $\partial^jz/\partial \tau^j$ is the result of the successive application $j$ times of $D^j w$ to the vector $\tau$, since $w \in P_k(\tilde{T})$, we can claim that
\begin{equation}
\label{dkzdtauhk}
\left\{
\begin{array}{l}
\displaystyle \frac{\mbox{d}^lz}{\mbox{d}\tau_h^l} = \displaystyle \sum_{j=1}^{\min[l,k]} a_{jl} \frac{\partial^j w}{\partial \tau^j} \\
\mbox{with} \\
\| a_{jl} \|_{2k-l,\infty,e_T} \leq C_{jl}(\Gamma) \mbox{ for } j=1,2,\ldots,\min[l.k].
\end{array}
\right.
\end{equation}
From the above observation together with \eqref{dkzdtauhk}, it becomes obvious that there exists a constant $C^a_{2k}(\Gamma)$ depending only on $\Gamma$ such that
\begin{equation}
\label{boundzw}
\| z \|_{2k,e_T} \leq C^a_{2k}(\Gamma) h_T^{1/2} \displaystyle \left(\| w \|^2_{0,\infty,\tilde{T}}+\sum_{j=1}^k \| D^j w \|^2_{0,\infty,\tilde{T}} \right)^{1/2} = C^a_{2k}(\Gamma) h_T^{1/2} \| w \|_{k,\infty,\tilde{T}}.
\end{equation}
Plugging \eqref{boundzw} into \eqref{chbarch4} and using an inverse inequality for $\tilde{T}$ given in \cite{ZAMM}, we come up with 
a constant $\bar{C}^a_H(\eta)$ depending only on $\eta$ and $\Gamma$, such that   
\begin{equation}
\label{chbarch5}
{\mathcal H}_T^{2k-1}(\bar{\mathcal G}^{a}_T(w) v) \leq \bar{C}^a_H(\eta) h_T^{k} \| w \|_{k,T}  \| v \|_{1/2,e_T}.
\end{equation}

To complete the proof we turn our attention to the term $\bar{\mathcal G}^{b}_T(w)$. \\
In this case we can take $\phi = \mbox{cos } \theta_T \bar{d}$ and $z = \bar{w_{\bf n}}$ where $w_{\bf n} = \partial w/\partial n$ 
is a function defined only on $\Gamma$. \\
As far as $\phi$ is concerned, there is nothing new, and actually, from \eqref{bhbarbh2ter} we infer the existence of a constant $\bar{C}_D(d)$ expressed in the same manner as $\bar{C}_H(\eta)$ by replacing $\eta$ with the trace of $d$ on $\Gamma$. Otherwise stated, it holds
\begin{equation}
\label{chbarch5bis}
{\mathcal H}^{2k-1}_T({\mathcal G}^b(w) v) \leq \bar{C}_D(d) h_T^{k+1/2} \| z \|_{2k,e_T} \| v \|_{1/2,e_T}.
\end{equation}
Here we resort again to mathematical induction, by making use of the property 
$$[\mbox{d}(\cdot)/\mbox{d}\tau_h](M) = [\mbox{cos }\theta_T \mbox{d}(\cdot)/\mbox{d}\tau](P).$$ 
More precisely, starting from the equation
$$[\mbox{d}^2 z/\mbox{d}\tau_h^2](M) = [-\mbox{sin }\theta_T \mbox{d}\theta_T/\mbox{d}\tau_h](M) [\mbox{d}w_{\bf n}/\mbox{d}\tau](P),$$
i.e.    
$$[\mbox{d}^2 z/\mbox{d}\tau_h^2](M) = [-y^{'}\kappa](M) [\mbox{d}w_{\bf n}/\mbox{d}\tau](P) + [1+(y^{'})^2]^{-1/2}(M)[\mbox{d}^2w_{\bf n}/\mbox{d}\tau^2](P).$$
we follow the same steps as in the case of $\bar{w}$ treated above, to infer the existence of constants $C^b_{jl}(\Gamma)$ depending only on $\Gamma$ for $1 \leq l \leq 2k$ and $j=1,\ldots,l$ such that 
$$|z|_{l,\infty,e_T} \leq \sum_{j=0}^{l} C^b_{jl}(\Gamma) | w_{\bf n} |_{j,\infty,\Gamma_T}.$$ 
This implies in turn the existence of a constant $C^b_{2k}(\Gamma)$ depending only on $\Gamma$, such that
\begin{equation}
\label{boundzwn}
\| z \|_{2k,e_T} \leq C^b_{2k}(\Gamma) h_T^{1/2} \displaystyle \left[\sum_{l=0}^{2k} | w_{\bf n} |^2_{l,\infty,\Gamma_T} \right]^{1/2}.
\end{equation}
Now noticing that $w_{\bf n} = \mbox{cos }\theta_T \partial w /\partial n_h + \mbox{sin }\theta_T \partial w /\partial \tau_h$, akin to  \eqref{dlzdtauhl} we can further establish the existence of constants $C^{\bf t}_{jl}(\Gamma)$ depending only on $\Gamma$, for $l=1,\ldots,2k$ and $j=1,\ldots,n$, with $n:=\min[l+1,k]$,such that 
\begin{equation}
\label{normaxdlwndtaul}
 | w_{\bf n} |^2_{l,\infty,\Gamma_T} \leq \displaystyle \sum_{j=1}^{n} C^{\bf t}_{jl} \| D^j w \|_{0,\infty,\tilde{T}}.
\end{equation}
Plugging \eqref{normaxdlwndtaul} into \eqref{boundzwn} and the resulting upper bound into \eqref{chbarch5bis}, we immediately come up with two mesh-independent constants $C^b_D(d)$ and $\bar{C}^b_D(d)$ such that
\begin{equation}
\label{chbarch6}
\left\{
\begin{array}{l}
{\mathcal H}_T^{2k-1}(\bar{\mathcal G}^{b}_T(w) v ) \leq  h_T^{k+1} \| v \|_{1/2,e_T} C^b_D(d) \| \nabla w \|_{k-1,\infty,\tilde{T}} \\  
\leq  h_T^{k} \| v \|_{1/2,e_T} \bar{C}_D^b(d) \| w \|_{k,T}.
\end{array}
\right.
\end{equation} 
Taking into account \eqref{chbarch6} and recalling \eqref{chbarch3} and \eqref{chbarch5}, for $\bar{C}_{C} = C_H(\eta) + C_{D}(\tilde{d})+ \bar{C}^a_H(\eta) + \bar{C}^b_D(d)$ we obtain 
\begin{equation}
\label{chbarch7}
\left\{
\begin{array}{l}
\bar{c}_T(w,v) - c_T(w,v) \leq \bar{C}_C h_T^{k} \| w \|_{k,T} \| v \|_{1/2,e_T}.
\end{array}
\right.
\end{equation}
Finally, summing up over ${\mathcal S}_h$ and using the Trace Theorem, we come up with a constant ${\mathcal C}_C$ depending on $\eta$, $\tilde{d}$ and $\Gamma$ but not on $h$, such that
\begin{equation}
\label{chbarch8}
|\bar{c}_h(w,v) -c_h(w,v)| \leq {\mathcal C}_C h^{k} \| w \|_{k,h} \| v \|_{1,h}.
\end{equation}
Combining \eqref{bhbarbh4} and \eqref{chbarch8} the result follows with ${\mathcal C}_A :=  {\mathcal C}_B + {\mathcal C}_C$.  
\QED \\

As a consequence, we can state the main result of this section, namely,
\begin{theorem}
\label{H1convergence}
Assume that the solution $u$ of \eqref{eq} belongs ot $H^{k+1}(\Omega)$. Let $\tilde{u}$ be an extension of $u$ in $H^{k+1}(\tilde{\Omega})$. Provided $h$ is sufficiently small and the quadrature formulae ${\mathcal J}^k_T$ and ${\mathcal I}^k_T$ integrate exactly polynomials of degree respectively less than or equal to $2k-2$ in $T$ and of degree less than or equal to $2k-1$ in $e_T$, the following error estimate holds for the solution $u^k_h$ of \eqref{vph} with a constant ${\mathcal C}_k$ independent of $h$:
\begin{equation}
\label{errestimukh}
\| \tilde{u} - u^k_k \|_{1,h} \leq {\mathcal C}_k h^k \| \tilde{u} \|_{k+1,\tilde{\Omega}}.
\end{equation} 
\end{theorem}
\prov
First take $w = \bar{u}^k_h$ in the estimate \eqref{ahbarah} and then use the second inequality of \eqref{estimatemk} to obtain an 
alternative estimate in terms of $\tilde{u}$. Finally, plugging the latter estimate into \eqref{errorbound0}, the result follows. \QED
 
\begin{remark}
The accuracy of the quadrature formulae ${\mathcal J}^k_T$ and ${\mathcal I}^k_T$, together with the regularity of $\Omega$ required for 
Theorem \ref{H1convergence} to hold may seem excessive. In this respect we first observe that, in practice, curved domains are often very smooth, such as disks, annuli and shapes alike. Moreover, we underline that $k$ equals $2$ or at most $3$ is overwhelmingly the usual choice of those practitioners, who refrain from using piecewise linear finite elements to approximate regular solutions, with the aim of gaining enhanced precision. More concretely, $2k-2$ and $2k-1$ equal just two and three for $k=2$ and hence a three-point Gauss integration formula in a triangle and a two-point Gauss integration formula in a segment are sufficient to ensure second order convergence. In case $k=3$ the use of a six-point quadrature formula in a triangle (cf. \cite{Dunavant}) and a three-point Gauss quadrature formula in a segment (cf. \cite{Quarteroni}) leads to optimal third order convergence. 
Nevertheless, in Section 6 we consider an alternative to cope with this problem by avoiding highly accurate quadrature formulae for larger values of $k$. Incidentally, we note that all the points for both aforementioned two-dimensional integration formulae lie inside the triangle, as required (see e.g. \cite{Dunavant}).  \QED
\end{remark}
 
\section{Numerical validation}
 
In this section we present the main results of the numerical experimentation carried out, in order to validate and check in different situations the theoretical studies for the method advocated in this work to solve second order elleptic equations with natural boundary conditions. 
Our numerical tests are confined to the case of quadratic finite elements. Moreover, for the sake of simplicity, we take $d \equiv 1$, $r \equiv 1$ and $\eta \equiv 0$ or $\eta \equiv 1$. Notice that, even with such a choice, in principle,  the use of numerical integration is necessary because of the trigonometric functions in the boundary corrective terms of the bilinear form $a_h$. \\
The comparisons are reported in two separate parts: In Subsection 5.1 we confront the performance of our new approach with the so-called do-nothing strategy, in which no corrective boundary terms are incorporated into the bilinear form, akin to problem \eqref{vp0h}. In Subsection 5.2 the performances in terms of accuracy of our method and the classical isoparametric technique (cf. \cite{BarrettElliott}) are compared.\\

More concretely, a certain number of test-problems with known exact solution satisfying different types of boundary conditions are solved. In all cases $\Omega$ is a domain contained in the ellipse $\Omega_1$ with semi-axes equal to $a$ and $b$, whose boundary is denoted by $\Gamma_1$. We consider both the case where $\Omega$ is $\Omega_1$ itself or the case where $\Omega$ is the hollow domain $\Omega_{hol} := \Omega_1 \setminus \overline{\Omega_{1/2}}$, $\Omega_{1/2}$ being the concentric ellipse with boundary $\Gamma_{1/2}$ and semi-axes $a/2$ and $b/2$ aligned with those of $\Omega_1$. Since we only considered test-problems with two axes of symmetry, we computed with quasi-uniform family of meshes for the quadrant ${\Omega}_1^{++}$ of $\Omega_1$ given by $x>0$ and $y>0$ generated by the procedure defined by an even integer $N$ described in \cite{CAMWA}. We recall that for $\Omega_1^{++}$ such meshes consist of $2N^2$ triangles. In case $\Omega$ is the hollow domain $\Omega_{hol}$ the final mesh is obtained by simply removing from the mesh of ${\Omega}_1^{++}$ the $N^2/2$ triangles fully contained in the corresponding quadrant of the ellipse $\bar{\Omega}_{1/2}$. \\
\indent Henceforth we use as mesh parameter the quantity $h_N:=1/N$. In all cases we take $N=2^n$ for $n$ ranging between $1$ and $4$.\\

For all the test-problems addressed in this section the results are summarized in tables with the following layout:\\
First we note that a given table displays results obtained with either the formulation \eqref{vph} or their counterparts determined by another method indicated in the table caption. Each table consists of three subsequent groups of rows containing absolute errors for increasing values of $N$ from left to right. These groups consist of two rows, with a sole exception with three rows, but for convenience they will be referred to in the sequel as "pairs" in all cases. The lowest row of each pair displays results obtained with the method indicated in the table caption for comparison with \eqref{vph}. On the upper rows in turn the corresponding results obtained with our method are supplied. On the rows of the two upper pairs we supply the approximation errors of $\nabla u$ and $u$ - or of $\tilde{u}$ and $\nabla \tilde{u}$ in the case of a hollow domain -, measured in the standard norm of $L^2(\Omega_h)$. On the rows of the third pair the maximum absolute error of the computed nodal values are supplied, thereby mimicking an error in the maximum norm by means of a discrete maximum semi-norm instead.\\ 

Throughout this section $\| \cdot \|_h$ and $| \cdot |_{1,h}$ denote the standard norm of $L^2(\Omega^{++}_h)$ and semi-norm of $H^1(\Omega^{++}_h)$, where $\Omega^{++}_h$ is the quarter sub-domain of $\Omega_h$ characterized by $x \geq 0$ and $y \geq 0$. In addition to this, $| \cdot |_{\infty,h}$ will represent the aforementioned discrete maximum semi-norm. For better readability, in the tables the numerical solution will be denoted by $u_N$ instead of $u^2_h$. \\

In all test cases reported in this section the underlying system of linear algebraic equations was solved by the classical Crout's method for band matrices. As for integrals along segments, whenever necessary, the following numerical quadrature formulae were employed on each boundary edge, depending on the case as specified hereafter: either the classical fourth order Simpson's rule referred to in the sequel as the $S$-formula or a three-point Gaussian formula of the sixth order referred to henceforth as the $G$-formula. On the other hand in all cases integrals in a straight-edged triangle are either determined exactly for polynomial integrands of degree less than or equal to two, or are approximated by a seven-point Gauss quadrature formula of the fifth order otherwise.  

\subsection{Comparisons with the do-nothing strategy}

Throughout this sub-section, instead of \eqref{eq}, we consider a reaction-advection-diffusion equation, by adding to the left hand side of the equation in $\Omega$ of \eqref{eq} the term $x \partial u/\partial x - y \partial u/\partial y$. \\
In the test-problems addressed below one-dimensional integrals were computed by using the $S$-formula \ for test-problems with Neumann boundary conditions, while the $G$-formula was rather used in the case of test-problems with Robin boundary conditions.\\

The counterpart of $u_N$ computed using on-site boundary conditions (i.e., without the corrections inherent to \eqref{vph}) referred to as the numerical solution for a do-nothing strategy, will be denoted by $u_{0N}$. \\
It should be noted that in all the test-problems addressed in this sub-section the tangential derivative of the exact solution vanishes identically on the boundary of the problem-definition domain. According to the a priori analysis given in Section 2, this condition is likely to ensure that $u_{0N}$ is a second order approximation of $u$ in the norm $\| \cdot \|_{1,h}$, which include a homogeneous Neumann problem such as \eqref{eq0}. Hence more accurate approximations $u_{0N}$ can be expected in this case. As seen below, in all cases indeed, the estimated order of $u_{0N}$ in this norm is two, akin to $u_N$, but the accuracy of the former is much poorer than the one of the latter, let alone the case of non zero tangential derivatives. Moreover the do-nothing approach remains suboptimal in terms of the two other measures under consideration.\\ 

\noindent \textbf{Test-problem 1 - consistency check for Neumann boundary conditions}\\

To begin with, we take $\eta=0$ and $\Omega$ equal to the unit disk (i.e. $\Omega = \Omega_1$ and $a=b=1$). 
Choosing $u = x^2 + y^2$ as an exact solution, we have $f \equiv -4 + 3x^2-y^2$ and $g \equiv 2$. Since $u$ is quadratic, our method is supposed to reproduce this solution exactly up to round-off errors. Indeed, the boundary integrals that come into play in the expression of $a_h$ for $\eta=0$ cancel out with the variational residual in this case, as one can infer from the arguments in Section 2 for problem \eqref{eq0}. On the other hand the do-nothing approach is not supposed to enjoy the same property. As a matter of fact, such effects are clearly observed in Table 1. As previously stated, the observed order of convergence of the latter strategy is two in the $| \cdot |_{1,h}$-semi-norm, but also in the $\| \cdot \|_h$-norm and in the $| \cdot |_{\infty,h}$-semi-norm, while third order convergence for both should rather show up. As seen from the results for Test-problem 2 hereafter, this confirms that the do-nothing strategy is indeed sub-optimal for $k=2$.\\

\begin{table*}[h!]
{\small 
\centering
\begin{tabular}{cccccc} &\\ [-.3cm]  
$h_N$ & $\longrightarrow$ & $1/2$ & $1/4$ & $1/8$ & $1/16$  
\tabularnewline &\\ [-.3cm] \hline &\\ [-.3cm]
$ | u_N - u |_{1,h}$ & $\longrightarrow$ & 0.2492264E-14 & 0.5707977E-14 & 0.1245960E-13 & 0.4093076E-13
\tabularnewline &\\ [-.3cm] \hline &\\ [-.3cm] 
$| u_{0N} - u|_{1,h}$ & $\longrightarrow$ & 0.23464732E-1 & 0.59076344E-2 & 0.14795511E-2 & 0.37005376E-3
\tabularnewline &\\ [-.3cm] \hline &\\ [-.3cm] 
\tabularnewline &\\ [-.3cm] \hline &\\ [-.3cm]
$\parallel u_N - u \parallel_h$ & $\longrightarrow$ & 0.1424145E-14 & 0.2398193E-13 & 0.1110170E-13 & 0.1636741E-12  
\tabularnewline &\\ [-.3cm] \hline &\\ [-.3cm]  
$\parallel u_{0N} - u \parallel_h$ & $\longrightarrow$ & 0.71108799E-1 & 0.17756383E-1 & 0.44378979E-2 & 0.11094013E-2 
\tabularnewline &\\ [-.3cm] \hline &\\ [-.3cm]
\tabularnewline &\\ [-.3cm] \hline &\\ [-.3cm] 
$| u_N - u |_{\infty,h}$ & $\longrightarrow$ & 0.2886580E-14 & 0.2986500E-13 & 0.1909584E-13 & 0.2112754E-12  
\tabularnewline &\\ [-.3cm] \hline &\\ [-.3cm] 
$| u_{0N} - u |_{\infty,h}$  & $\longrightarrow$ & 0.93503269E-1 & 0.23091994E-1 & 0.57556584E-2 & 0.14378358E-2
\tabularnewline &\\ [-.3cm] \hline &\\ [-.3cm] 
\end{tabular}
\caption{Errors for Test-problem 1 solved by the method \eqref{vph} and the do-nothing approach with $k=2$} 
}
\label{table1}
\end{table*}

\noindent \textbf{Test-problem 2 - consistency check for Robin boundary conditions} \\

Next we consider the same data as in Test-problem 1, except for $\eta$ and $g$. The former now equals one, while $g = 2 + u_{|\Gamma} \equiv 3$. In spite of the fact that $u$ is quadratic, our method here fails to yield this solution exactly (up to round-off errors). This is because the term $\int_{e_T} \eta \mbox{ cos }\theta_T \bar{w} v$ of the element matrix for a triangle $T \in {\mathcal S}_h$ cannot be computed exactly at all by any quadrature formula of the form ${\mathcal I}^k_T$ (cf. \eqref{ITk}). Nevertheless, while the do-nothing approach behaves basically in the same sub-optimal manner as in the previous test-problem, the observed orders of convergence of our method are optimal, namely, equal to three in both the $\| \cdot \|_h$-norm and  the $| \cdot |_{\infty,h}$-semi-norm and roughly equal to 5/2 in the $| \cdot |_{1,h}$-semi-norm, as one can infer from Table 2. We have no explanation for the latter order half a point beyond expectancy, except perhaps for the particularity of this test-case.\\

\begin{table*}[h!]
{\small 
\centering
\begin{tabular}{cccccc} &\\ [-.3cm]  
$h_N$ & $\longrightarrow$ & $1/2$ & $1/4$ & $1/8$ & $1/16$  
\tabularnewline &\\ [-.3cm] \hline &\\ [-.3cm]
$| u_N - u |_{1,h}$ & $\longrightarrow$ & 0.26559027E-2 & 0.49169500E-3 & 0.88849622E-4 & 0.15872499E-4
\tabularnewline &\\ [-.3cm] \hline &\\ [-.3cm] 
$| u_{0N} - u |_{1,h}$ & $\longrightarrow$ & 0.12178188E-1 & 0.30486276E-2 & 0.75512677E-3 & 0.18739613E-3
\tabularnewline &\\ [-.3cm] \hline &\\ [-.3cm] 
\tabularnewline &\\ [-.3cm] \hline &\\ [-.3cm]
$\parallel u_N - u \parallel_h$ & $\longrightarrow$ & 0.39908099E-3 & 0.50529798E-4 & 0.62954104E-5 & 0.78335946E-6  
\tabularnewline &\\ [-.3cm] \hline &\\ [-.3cm]  
$\parallel u_{0N} - u \parallel_h$ & $\longrightarrow$ & 0.35172026E-1 & 0.88651467E-2 & 0.22207750E-2 & 0.55547324E-3 
\tabularnewline &\\ [-.3cm] \hline &\\ [-.3cm]
\tabularnewline &\\ [-.3cm] \hline &\\ [-.3cm] 
 $| u_N - u |_{\infty,h}$ & $\longrightarrow$ & 0.72665131E-3 & 0.10073351E-3 & 0.13067108E-4 & 0.16640715E-5  
\tabularnewline &\\ [-.3cm] \hline &\\ [-.3cm] 
$| u_{0N} - u |_{\infty,h}$  & $\longrightarrow$ & 0.45850174E-1 & 0.11461358E-1 & 0.28615773E-2 & 0.71464685E-3
\tabularnewline &\\ [-.3cm] \hline &\\ [-.3cm] 
\end{tabular}
\caption{Errors for Test-problem 2 solved by the method \eqref{vph} and the do-nothing approach with $k=2$}
}
\label{table2}
\end{table*}

\noindent \textbf{Test-problem 3 - order and accuracy check for a convex domain} \\

The aim of this test-problem and the next one is to check the estimated orders of convergence of our method and the do-nothing strategy by solving a Robin problem with a non polynomial exact solution. Additionally we compare the accuracy provided by both approaches. Still taking $\eta=1$ now $\Omega$ is the ellipse $\Omega_1$ with $a=1/2$ and $b=1$. 
Setting $\rho=\sqrt{x^2/a^2+y^2/b^2}$ we take the exact solution $u = \rho^3/3$ and determine $f$ and $g$ accordingly. Here, although $u_{|\Gamma} \equiv 1/3$, $g$ is not constant on $\Gamma$ since $(\partial u/\partial n)_{|\Gamma} = \sqrt{x^2/a^4+y^2/b^2}$. The errors for this test-problem are supplied in Table 3, from which we infer the order of convergence of ca. two in the $| \cdot |_{1,h}$-semi-norm for both numerical techniques being compared. The observed order of convergence for our method in turn equals three both in the $\| \cdot \|_h$-norm and in the $| \cdot |_{\infty,h}$-semi-norm, while it remains about two for the do-nothing approach. Furthermore, akin to the two preceding test-problem, our method turns out to be significantly more accurate than the latter.\\
          
\begin{table*}[h!]
{\small 
\centering
\begin{tabular}{cccccc} &\\ [-.3cm]  
$h_N$ & $\longrightarrow$ & $1/2$ & $1/4$ & $1/8$ & $1/16$  
\tabularnewline &\\ [-.3cm] \hline &\\ [-.3cm]
$ | u_N - u |_{1,h}$ & $\longrightarrow$ & 0.13579836E-1 & 0.31352173E-2 & 0.77561168E-3 & 0.19128506E-3
\tabularnewline &\\ [-.3cm] \hline &\\ [-.3cm] 
$| u_{0N} - u|_{1,h}$ & $\longrightarrow$ & 0.28251034E-1 & 0.80090351E-2 & 0.20785941E-2 & 0.52505776E-3
\tabularnewline &\\ [-.3cm] \hline &\\ [-.3cm] 
\tabularnewline &\\ [-.3cm] \hline &\\ [-.3cm]
$\parallel u_N - u \parallel_h$ & $\longrightarrow$ & 0.51615327E-2  & 0.35252677E-3 & 0.34528329E-4 & 0.40299932E-5  
\tabularnewline &\\ [-.3cm] \hline &\\ [-.3cm]  
$\parallel u_{0N} - u \parallel_h$ & $\longrightarrow$ &  0.17927294E-1 & 0.50097716E-2 & 0.12871556E-02 & 0.32387353E-3 
\tabularnewline &\\ [-.3cm] \hline &\\ [-.3cm]
\tabularnewline &\\ [-.3cm] \hline &\\ [-.3cm] 
$| u_N - u |_{\infty,h}$ & $\longrightarrow$ & 0.24067504E-1 & 0.22750493E-2 &  0.22424154E-3 & 0.25823392E-4  
\tabularnewline &\\ [-.3cm] \hline &\\ [-.3cm] 
$| u_{0N} - u |_{\infty,h}$  & $\longrightarrow$ &   0.68104252E-1   & 0.18825715E-1 & 0.47541530E-2 & 0.11868257E-2
\tabularnewline &\\ [-.3cm] \hline &\\ [-.3cm] 
\end{tabular}
\caption{Errors for Test-problem 3 solved by the method \eqref{vph} and the do-nothing approach with $k=2$} 
}
\label{table3}
\end{table*}

\noindent \textbf{Test-problem 4 - order and accuracy check for a non convex domain}\\

Still taking $\eta=1$ now $\Omega$ is the annulus $\Omega_{hol}$ with $a=b=1$. 
Setting $\rho_I=1/2$, $\rho_E=1$ and $\rho=\sqrt{x^2+y^2}$ we consider $u = \rho^3/(3 \rho_I \rho_E) - \rho^2(\rho_I^{-1}+\rho_E^{-1})/2 + \rho$ as an exact solution and determine $f$ and $g$ accordingly. Incidentally $g$ is constant on each boundary portion given by $\rho = \rho_E$ and $\rho = \rho_I$. The errors for this test-problem are supplied in Table 4. From these data we infer orders of convergence close to two in the $| \cdot |_{1,h}$-semi-norm for both approaches. On the other hand, this estimated order increases to ca. three for our method and remains about two for the do-nothing approach, both in the $\| \cdot \|_h$-norm and in the $| \cdot |_{\infty,h}$-semi-norm, akin to Test-problem 3. Moreover our method is much more accurate than the latter all the way. \\
In Table 4 $\tilde{u}$ stands for the function equal to $u$ in $\Omega$, such that $\tilde{u}(x,y)$ is expressed in the same way as $u(x,y)$. 
in $\Omega_h \setminus \Omega$.  \\
                                     
\begin{table*}[h!]
{\small 
\centering
\begin{tabular}{cccccc} &\\ [-.3cm]  
$h_N$ & $\longrightarrow$ & $1/2$ & $1/4$ & $1/8$ & $1/16$  
\tabularnewline &\\ [-.3cm] \hline &\\ [-.3cm]
$ | u_N - \tilde{u} |_{1,h}$ & $\longrightarrow$ & 0.10529838E-1 & 0.33293060E-2 & 0.90693270E-3 &  0.23429414E-3
\tabularnewline &\\ [-.3cm] \hline &\\ [-.3cm] 
$| u_{0N} - \tilde{u}|_{1,h}$ & $\longrightarrow$ & 0.87486680E-1 & 0.22599845E-1 & 0.56444817E-2 &  0.14040942E-2
\tabularnewline &\\ [-.3cm] \hline &\\ [-.3cm] 
\tabularnewline &\\ [-.3cm] \hline &\\ [-.3cm]
$\parallel u_N - \tilde{u} \parallel_h$ & $\longrightarrow$ & 0.11286345E-1 & 0.15298096E-2 & 0.20047702E-3 & 0.25763950E-4
\tabularnewline &\\ [-.3cm] \hline &\\ [-.3cm]  
$\parallel u_{0N} - \tilde{u} \parallel_h$ & $\longrightarrow$ & 0.87410070E-1 & 0.22987131E-1 & 0.58163990E-2 & 0.14583711E-2
\tabularnewline &\\ [-.3cm] \hline &\\ [-.3cm]
\tabularnewline &\\ [-.3cm] \hline &\\ [-.3cm] 
$| u_N - \tilde{u} |_{\infty,h}$ & $\longrightarrow$ & 0.17669426E-01 & 0.22958502E-2 &   0.29871585E-3 & 0.38406792E-4  
\tabularnewline &\\ [-.3cm] \hline &\\ [-.3cm] 
$| u_{0N} - \tilde{u} |_{\infty,h}$  & $\longrightarrow$ & 0.15082914E+0 & 0.38490394E-1 &  0.96152869E-2  &  0.23960461E-2
\tabularnewline &\\ [-.3cm] \hline &\\ [-.3cm] 
\end{tabular}
\caption{Errors for Test-problem 4 solved by the method \eqref{vph} and the do-nothing approach with $k=2$} 
}
\label{table4}
\end{table*}

\subsection{Comparisons with the isoparametric technique}

Throughout this sub-section the problem to be solved is \eqref{eq}. In this manner, the comparisons with the isoparametric method are strictly restricted to the analytical framework considered both in \cite{BarrettElliott} and in the present work. \\
Numerical integration in an edge of the master triangle $\hat{T}$ is performed by means of the classical two-point Gauss quadrature formula of the fourth order for the isoparametric computations. On the other hand, the $G$-formula is employed for one-dimensional integrals in the case of our method. However, in the latter case, for comparative purposes, the computations with our method for solving Test-problem 7 were also carried out by using the $S$-formula as well, for the same purpose. \\
Throughout this subsection the domain $\Omega$ is the ellipse $\Omega_1$ with semi-axes $a=1/2$ and $b=1$. \\

In the tables that follow the numerical solution obtained by the isoparametric method is denoted by $\tilde{u}_N$ for increasing values of $N$.\\  

\noindent \textbf{Test-problem 5 - consistency check for Neumann boundary conditions}\\

Choosing $u = 4x^2 + y^2$ as an exact solution, we have $f = 4x^2 + y^2 - 10$ and $g = 2 \sqrt{16x^2 + y^2}$. Since $u$ is quadratic, our method is supposed to reproduce this solution exactly up to round-off errors. Indeed, the boundary integrals that come into play in the expression of $a_h$ for $\eta=0$ cancel out with the variational residual in this case, as one can infer from the arguments in Section 2 for problem \eqref{eq0}. On the other hand, it is noticeable that the isoparametric approach does not enjoy such a property, since it can only represent exactly a quadratic function in a current triangle $T$, if the mapping from the master triangle onto it is affine. As a matter of fact, such effects are clearly observed in Table 5. As for isoparametrics the estimated order of convergence is roughly $5/2$ in the $| \cdot |_{1,h}$-semi-norm, quite unexpectedly for $k=2$. Super-convergence in the $\| \cdot \|_h$-norm, more precisely of order close to $7/2$, is also observed. However, the rate of convergence in the $| \cdot |_{\infty,h}$-semi-norm is sub-optimal, since it is ca. two. Moreover, the accuracy itself in terms of the latter measure is rather poor. The same issues come out in the following test-problems, for which a conjectural explanation will be given afterward.\\ 

\begin{table*}[h!]
{\small 
\centering
\begin{tabular}{cccccc} &\\ [-.3cm]  
$h_N$ & $\longrightarrow$ & $1/2$ & $1/4$ & $1/8$ & $1/16$  
\tabularnewline &\\ [-.3cm] \hline &\\ [-.3cm]
$ | u_N - u |_{1,h}$ & $\longrightarrow$ & 0.2729988E-14 & 0.4184940E-14 & 0.1561916E-13 & 0.8922446E-13
\tabularnewline &\\ [-.3cm] \hline &\\ [-.3cm] 
$| \tilde{u}_N - u|_{1,h}$ & $\longrightarrow$ & 0.25294001E-1 & 0.48526314E-2 & 0.87514425E-3 & 0.15603152E-3
\tabularnewline &\\ [-.3cm] \hline &\\ [-.3cm] 
\tabularnewline &\\ [-.3cm] \hline &\\ [-.3cm]
$\parallel u_N - u \parallel_h$ & $\longrightarrow$ & 0.2281721E-14 & 0.1741791E-14 & 0.9102329E-13 & 0.4946333E-12  
\tabularnewline &\\ [-.3cm] \hline &\\ [-.3cm]  
$\parallel  \tilde{u}_N - u \parallel_h$ & $\longrightarrow$ & 0.14213778E-2 & 0.12733116E-3 &  0.10879178E-4 & 0.95442362E-6 
\tabularnewline &\\ [-.3cm] \hline &\\ [-.3cm]
\tabularnewline &\\ [-.3cm] \hline &\\ [-.3cm] 
$| u_N - u |_{\infty,h}$ & $\longrightarrow$ & 0.4884981E-14 & 0.4107825E-14 & 0.1587619E-12 & 0.8637535E-12  
\tabularnewline &\\ [-.3cm] \hline &\\ [-.3cm] 
$| \tilde{u}_N - u |_{\infty,h}$  & $\longrightarrow$ & 0.11704065E+0 & 0.35893397E-1 & 0.94800122E-2 & 0.24014197E-2
\tabularnewline &\\ [-.3cm] \hline &\\ [-.3cm] 
\end{tabular}
\caption{Errors for Test-problem 5 solved by the method \eqref{vph} and isoparametrics with $k=2$} 
}
\label{table5}
\end{table*}

\noindent \textbf{Test-problem 6 - consistency check for Robin boundary conditions}\\

Next we consider the same data as in Test-problem 5, except for $\eta$ and $g$. The former now equals one, thereby yielding $g = u +2 \sqrt{16x^2 + y^2}$, either on $\Gamma$ or on the piecewise parabolic boundary $\tilde{\Gamma}_h$ approximating $\Gamma$ in the way inherent to the isoparametric method. For the same reason as for Test-problem 2, in spite of the fact that $u$ is quadratic, in this case too our method is unable to reproduce this solution exactly (up to round-off errors). Here, the observed order of convergence of both our method and isoparametrics is roughly two, while both show the same accuracy as far as the $| \cdot |_{1,h}$-seminorm is concerned, with a slight advantage of the former over the latter. On the other hand, the behavior of both methods differs significantly in terms of the $\| \cdot \|_h$-norm and the $| \cdot |_{\infty,h}$-semi-norm. As for the former the isoparametric approach is some orders of magnitude more accurate, although both methods are roughly of the third order. In contrast, from the point of view of maximum errors, our method is much more accurate. As a matter of fact, the observed orders of convergence in terms of the $| \cdot |_{\infty,h}$-semi-norm seems to tend slowly to three for our approach and sub-optimally to two for the isoparametric technique, as one infers from Table 6.\\

\begin{table*}[h!]
{\small 
\centering
\begin{tabular}{cccccc} &\\ [-.3cm]  
$h_N$ & $\longrightarrow$ & $1/2$ & $1/4$ & $1/8$ & $1/16$  
\tabularnewline &\\ [-.3cm] \hline &\\ [-.3cm]
$| u_N - u |_{1,h}$ & $\longrightarrow$ & 0.14863869E-1 & 0.30569293E-2 & 0.57828450E-3 & 0.10524467E-3
\tabularnewline &\\ [-.3cm] \hline &\\ [-.3cm] 
$| \tilde{u}_N - u |_{1,h}$ & $\longrightarrow$ & 0.16986880E-1 & 0.31989699E-2 & 0.58082760E-3 & 0.10430537E-3
\tabularnewline &\\ [-.3cm] \hline &\\ [-.3cm] 
\tabularnewline &\\ [-.3cm] \hline &\\ [-.3cm]
$\parallel u_N - u \parallel_h$ & $\longrightarrow$ & 0.16275044E-2 & 0.28902417E-3 & 0.47903737E-4 & 0.63583308E-5  
\tabularnewline &\\ [-.3cm] \hline &\\ [-.3cm]  
$\parallel  \tilde{u}_N - u \parallel_h$ & $\longrightarrow$ & 0.74702192E-3 & 0.72138092E-4 & 0.67781245E-5 & 0.62199865E-6 
\tabularnewline &\\ [-.3cm] \hline &\\ [-.3cm]
\tabularnewline &\\ [-.3cm] \hline &\\ [-.3cm] 
$| u_N - u |_{\infty,h}$ & $\longrightarrow$ & 0.23987536E-1 & 0.16344968E-2 & 0.33411802E-3 & 0.46717335E-4  
\tabularnewline &\\ [-.3cm] \hline &\\ [-.3cm] 
$| \tilde{u}_N - u |_{\infty,h}$  & $\longrightarrow$ & 0.10259756E+0 & 0.34089002E-1 & 0.93091702E-2 & 0.23854162E-2
\tabularnewline &\\ [-.3cm] \hline &\\ [-.3cm] 
\end{tabular}
\caption{Errors for Test-problem 6 solved by the method \eqref{vph} and isoparametrics with $k=2$}
}
\label{table6}
\end{table*}

\noindent \textbf{Test-problem 7 - order and accuracy check for homogeneous Neumann boundary conditions}\\

In solving this test-problem we estimate the orders of convergence of our method and the isoparametric technique in the approximation of a non polynomial exact solution, while comparing the accuracy of both approaches. Moreover here we also confront two different one-dimensional numerical quadrature formulae to compute the terms of the element matrices for triangles in ${\mathcal S}_h$ accounting for the integrals along $\Gamma_h$ inherent to our method. The solutions obtained by using the $G$-formula are still denoted by $u_N$, while those computed by employing Simpson's rule are denoted by $u_N^S$. \\
Taking $\eta=0$ and setting $\rho=\sqrt{4x^2+y^2}$ here the manufactured exact solution is given by $u = \rho^2- 2\rho^3/3$, $f$ being  determined accordingly. Since $g \equiv 0$ on $\Gamma$, this value is extended to the boundary's immediate neighborhood in order to implement the isoparametric solution method. The errors for this test-problem are supplied in Table 7, from which we infer the order of convergence of ca. two in the $| \cdot |_{1,h}$-semi-norm for the three numerical techniques being compared. As far as accuracy is concerned, our method shows up superior, and even more with the use of the $S$-formula. 
In the case of the $\| \cdot \|_h$-norm in turn, again the observed order of convergence is about three for both methods being compared in all versions, while the isoparametric technique turns out to be significantly more accurate. Furthermore, in contrast to the two previous test-problems, the maximum errors measured in the $| \cdot |_{\infty,h}$-semi-norm are also smaller for the isoparametric approach. Actually here the order of the latter method in terms of this semi-norm increases to an optimal three, thereby equaling our method for both $u_N$ and $u_N^S$. Such a different behavior could be due to the fact that $\eta=0$ and $g \equiv 0$ in the present case, a conjecture reinforced by the results for the next test-problem. \\
          
\begin{table*}[h!]
{\small 
\centering
\begin{tabular}{cccccc} &\\ [-.3cm]  
$h_N$ & $\longrightarrow$ & $1/2$ & $1/4$ & $1/8$ & $1/16$  
\tabularnewline &\\ [-.3cm] \hline &\\ [-.3cm]
$ | u_N - u |_{1,h}$ & $\longrightarrow$ & 0.23864345E-1 & 0.56306694E-2 & 0.14328501E-2 & 0.36515127E-3
\tabularnewline &\\ [-.3cm] \hline &\\ [-.3cm]
$ | u_N^S - u |_{1,h}$ & $\longrightarrow$ & 0.33876440E-1 & 0.91820842E-2 & 0.23858551E-2 & 0.60877768E-3
\tabularnewline &\\ [-.3cm] \hline &\\ [-.3cm]  
$| \tilde{u}_{N} - u|_{1,h}$ & $\longrightarrow$ & 0.44748258E-1 & 0.12514075E-1 & 0.32754274E-2 & 0.83446142E-3
\tabularnewline &\\ [-.3cm] \hline &\\ [-.3cm] 
\tabularnewline &\\ [-.3cm] \hline &\\ [-.3cm]
$\parallel u_N - u \parallel_h$ & $\longrightarrow$ & 0.36893575E-1  & 0.46456671E-2 & 0.57989566E-3 & 0.72995933E-4  
\tabularnewline &\\ [-.3cm] \hline &\\ [-.3cm]  
$\parallel u_N^S - u \parallel_h$ & $\longrightarrow$ & 0.22924642E-1 & 0.30466850E-2 & 0.39519880E-3 & 0.50395917E-4  
\tabularnewline &\\ [-.3cm] \hline &\\ [-.3cm]  
$\parallel \tilde{u}_{N} - u \parallel_h$ & $\longrightarrow$ & 0.19714167E-2 & 0.27536829E-3 & 0.36247891E-4 & 0.46144765E-5 
\tabularnewline &\\ [-.3cm] \hline &\\ [-.3cm]
\tabularnewline &\\ [-.3cm] \hline &\\ [-.3cm] 
 $| u_N - u |_{\infty,h}$ & $\longrightarrow$ & 0.77933345E-1 & 0.99002832E-2 & 0.12402301E-2 & 0.15705859E-3  
\tabularnewline &\\ [-.3cm] \hline &\\ [-.3cm] 
$| u_N^S - u |_{\infty,h}$ & $\longrightarrow$ & 0.30112097E-1 & 0.40448281E-2 & 0.52908458E-3 & 0.69288886E-4  
\tabularnewline &\\ [-.3cm] \hline &\\ [-.3cm] 
$| \tilde{u}_{N} - u |_{\infty,h}$  & $\longrightarrow$ & 0.70182481E-2 & 0.11226893E-2 & 0.15730104E-3 & 0.20655585E-4
\tabularnewline &\\ [-.3cm] \hline &\\ [-.3cm] 
\end{tabular}
\caption{Errors for Test-problem 7 solved by the method \eqref{vph} and isoparametrics with $k=2$} 
}
\label{table7}
\end{table*}
\noindent \textbf{Test-problem 8 - order and accuracy check for inhomogeneous Robin boundary conditions}\\

In the last experiment reported in this section, we estimate the orders of convergence of our method and the isoparametric technique in the approximation of a non polynomial exact solution, while comparing the accuracy of both approaches. Taking $\eta=1$ and setting $\rho=\sqrt{x^2+y^2}$ we take the exact solution $u = \rho^2- 2\rho^3/3$ and determine accordingly $f = -4 + 6\rho + \rho^2 -2\rho^3/3$ and $g = (4x^2+y^2)(1-\rho)/\sqrt{16x^2+y^2}+u$. The errors for this test-problem are supplied in Table 8, from which we infer the order of convergence of ca. two in the $| \cdot |_{1,h}$-semi-norm for both numerical techniques being compared, with slightly improved accuracy for our method. 
In the case of the $\| \cdot \|_h$-norm in turn, the observed order of convergence is about three for both methods, but here again the isoparametric technique turns out to be significantly more accurate. In contrast, the maximum errors measured in the $| \cdot |_{\infty,h}$-semi-norm are much larger for the isoparametric method. Incidentally, in this respect, even more clearly than for Test-problem 6, the estimated order of convergence of our method is roughly three, while it is only ca. two for isoparametrics. Taking into account these  observations, such a discrepancy could be explained by singularities of the piecewise rational shape-functions inherent to the isoparametric approach, as compared to the more gentle behavior of the piecewise quadratics employed by our method, especially in the computation of one-dimensional integrals.\\
          
\begin{table*}[h!]
{\small 
\centering
\begin{tabular}{cccccc} &\\ [-.3cm]  
$h_N$ & $\longrightarrow$ & $1/2$ & $1/4$ & $1/8$ & $1/16$  
\tabularnewline &\\ [-.3cm] \hline &\\ [-.3cm]
$ | u_N - u |_{1,h}$ & $\longrightarrow$ & 0.19935758E-1 & 0.53764349E-2 & 0.14176285E-2 & 0.36420729E-3
\tabularnewline &\\ [-.3cm] \hline &\\ [-.3cm] 
$| \tilde{u}_N - u|_{1,h}$ & $\longrightarrow$ & 0.22407214E-1 & 0.59934938E-2 & 0.15136095E-2 & 0.37767148E-3
\tabularnewline &\\ [-.3cm] \hline &\\ [-.3cm] 
\tabularnewline &\\ [-.3cm] \hline &\\ [-.3cm]
$\parallel u_N - u \parallel_h$ & $\longrightarrow$ & 0.69836150E-2 & 0.87032627E-3 & 0.10798977E-3 & 0.13531435E-4  
\tabularnewline &\\ [-.3cm] \hline &\\ [-.3cm]  
$\parallel \tilde{u}_N - u \parallel_h$ & $\longrightarrow$ & 0.10705974E-2 & 0.14649138E-3 & 0.18424592E-4 & 0.22837888E-5 
\tabularnewline &\\ [-.3cm] \hline &\\ [-.3cm]
\tabularnewline &\\ [-.3cm] \hline &\\ [-.3cm] 
 $| u_N - u |_{\infty,h}$ & $\longrightarrow$ & 0.19711747E-1 & 0.26755367E-2 & 0.34243483E-3 & 0.44022901E-4  
\tabularnewline &\\ [-.3cm] \hline &\\ [-.3cm] 
$| \tilde{u}_N - u |_{\infty,h}$  & $\longrightarrow$ & 0.99880700E-2 & 0.21099200E-2 & 0.53693035E-3 & 0.13478228E-3
\tabularnewline &\\ [-.3cm] \hline &\\ [-.3cm] 
\end{tabular}
\caption{Errors for Test-problem 8 solved by the method \eqref{vph} and isoparametrics with $k=2$} 
}
\label{table8}
\end{table*}

\section{Final comments}

To conclude we make a few relevant comments on this work.

\begin{enumerate}
\item
\textbf{Main merits and demerits -} Improved accuracy can be expected from our method with respect to more classical approaches, such as isoparametric finite elements, since its key principle is the fulfillment of the right natural boundary conditions at the exact location where they are prescribed. 
While in this work such a conjecture turned out to be true in many important aspects, it remains to be 
confirmed through additional experimentation. Nevertheless, there is at least one a priori main advantage of our approach as compared to those found in the literature (see e.g. \cite{BarrettElliott} and references therein): Since the basis functions for the underlying approximation space are polynomials in each triangle, it is perfectly possible to do without numerical quadrature in order to compute element matrices and vectors, as long as the functions standing for problem data are replaced by their polynomial interpolates of a suitable degree in each element or on a boundary edge thereof if applicable. On the other hand, akin to the case of our method to handle Dirichlet boundary conditions prescribed on smooth bondaries (cf. \cite{ZAMM,IMAJNA} among other papers listed in the bibliography), the solution technique advocated here gives rise to a non symmetric matrix, even when the differential operator for the problem to solve is self-adjoint. Nevertheless, in the present case, it is possible to use as well the easy-to-implement and fast converging iterative procedure proposed in \cite{JCAMbis}, in which a linear system is solved at every step with the same symmetric positive definite matrix associated with the do-nothing approach, provided the underlying differential operator is of the aforementioned type.
\item
\textbf{A word on conditioning -} The previous remark gives a hint on the matrix condition number for methods based on the off-site enforcement of boundary conditions, such as the one studied here. Indeed, presuming that the aforementioned iterative procedure is fast converging in the present case too, it can be asserted that the conditioning of our method is just a small perturbation of the one applying to the do-nothing approach. It turns out that this is known to be fine, especially for self-adjoint coercive problems.  
\item 
\textbf{Regularity considerations -} The regularity required for the problem data and geometry in our reliability analysis summarized in the statement of Theorem \ref{H1convergence} can be considerably weakened, if instead of the quadrature formulae ${\mathcal J}_T^k$and ${\mathcal I}_T^k$ we use standard Lagrange interpolation at suitable points inside $T$ and $e_T$ to approximate the functions $d$, $r$, $\bar{g}$ and $\bar{w}$ in the definition of the approximate problem \eqref{vph}. In this case the error analysis would require only $C^k$-regularity from these functions and $\Omega$ itself, to yield the same qualitative results. However, we did not adopt such a strategy here, because we chose to stick as much as possible to the functional framework already exploited by the authors quoted in the bibliography, who had studied the same kind of problem prior to us. This allowed us to shorten to a non negligible extent the analysis carried out in this article.
\item 
\textbf{Accuracy vs. complexity -} 
A handful of representative tests with quadratic Lagrange finite elements were reported in Section 5. 
From the author's point of view, they illustrated how competitive the alternative advocated here to tackle the problem at hand is, as compared to more classical techniques studied in work on the topic prior to this one, such as isoparametric elements. In spite of this, it is true that the aforementioned tests raised a few issues, which cannot be fully explained by the a priori reliability analysis carried out in Sections 2 and 4. More particularly, we could underline the super-convergence observed in some cases for the approach doing without boundary corrective terms, referred to here as the do-nothing strategy. Notice however that, although, the latter technique is tempting due to its simplicity, the implementation of the method proposed here is just a little more complex. Moreover, the former appeared to be orders of magnitude less accurate than the latter for quadratic elements, let alone higher order methods. 
\item
\textbf{Perspectives for future work -} The material presented in this article encourages the author to push further formal and numerical studies on the off-site enforcement of non essential boundary conditions, in a natural and straightforward way at a time in the finite-element environment and beyond. More particularly he means their extension to the three-dimensional case and also to fourth order boundary value problems. Complements inherent to the second order  problems in the two-dimensional case studied in this work itself could also be addressed, such as computational tests with values of $k$ greater than two.
\end{enumerate}

\end{document}